\documentclass[review]{elsarticle}

\UseRawInputEncoding
\usepackage{lineno,hyperref}
\usepackage{graphicx}
\usepackage{subfigure}
\usepackage{amsmath}
\usepackage{amssymb}
\usepackage{txfonts}
\usepackage{pifont}
\usepackage{fleqn}
\usepackage{natbib}
\usepackage{geometry}
\usepackage{stmaryrd}
\usepackage{mathrsfs}
\usepackage{setspace}
\usepackage{enumerate}
\usepackage{threeparttable}

\usepackage[font=small,labelfont=bf,labelsep=none]{caption}  
\modulolinenumbers[5]

\newtheorem{thm}{Theorem}

\newtheorem{asmp}{Assumption}
\newdefinition{rmk}{Remark}
\newdefinition{exm}{Example}
\newtheorem{df}{Definition}
\newproof{proof}{\bf{Proof}}

\journal{Aerospace Science and Technology}

\biboptions{sort&compress}

\begin{document}
\captionsetup[figure]{labelfont={bf},labelformat={default},labelsep=period,name={Fig.}}
\begin{frontmatter}

\title{Penetration trajectory optimization for the hypersonic gliding vehicle encountering two interceptors\tnoteref{mytitlenote}}
\tnotetext[mytitlenote]{This work was supported by the Science and Technology Innovation 2030-Key Project of "New Generation Artificial Intelligence" under Grant 2020AAA0108200, the National Natural Science Foundation of China under Grants 61873011, 61973013, 61922008, 61803014, 62103023, and 62103016, the Defense Industrial Technology Development Program under Grant JCKY2019601C106, the Innovation Zone Project under Grant 18-163-00-TS-001-001-34, the Foundation Strengthening Program Technology Field Fund under Grant 2019-JCJQ-JJ-243, the Fund from Key Laboratory of Dependable Service Computing in Cyber Physical Society under Grant CPSDSC202001, the China National Postdoctoral Program for Innovative Talents under Grant BX20200034, and Project funded by China Postdoctoral Science Foundation under Grant 2020M680297.}


\author[mymainaddress]{Zhipeng Shen}

\author[mymainaddress]{Jianglong Yu\corref{mycorrespondingauthor}}
\cortext[mycorrespondingauthor]{Corresponding author}
\ead{sdjxyjl@buaa.edu.cn}

\author[b,mymainaddress,c]{Xiwang Dong}
\author[b]{Yongzhao Hua}
\author[mymainaddress,b,c]{Zhang Ren}

\address[mymainaddress]{School of Automation Science and Electrical Engineering, Science and Technology on Aircraft Control Laboratory, Beihang University, Beijing, 100191, P.R. China}
\address[b]{Institute of Artificial Intelligence, Beihang University, Beijing 100191, China}
\address[c]{Beijing Advanced Innovation Center for Big Data and Brain Computing, Beihang University, Beijing 100191, China}

\begin{abstract}
The penetration trajectory optimization problem for the hypersonic gliding vehicle (HGV) encountering two interceptors is investigated. The HGV penetration trajectory optimization problem considering the terminal target area is formulated as a nonconvex optimal control problem. The nonconvex optimal control problem is transformed into a second-order cone programming (SOCP) problem, which can be solved by state-of-the-art interior-point methods. In addition, a penetration strategy that only requires the initial line-of-sight angle information of the interceptors is proposed. The convergent trajectory obtained by the proposed method allows the HGV to evade two interceptors and reach the target area successfully. Furthermore, a successive SOCP method with a variable trust region is presented, which is critical to balancing the trade-off between time consumption and optimality. Finally, the effectiveness and performance of the proposed method are verified by numerical simulations.
\end{abstract}

\begin{keyword}
\texttt Hypersonic gliding vehicle\sep Penetration trajectory optimization\sep Second-order cone programming \sep Variable trust region \sep Interceptors
\end{keyword}

\end{frontmatter}


\section{Introduction}\label{section1}
The hypersonic gliding vehicle (HGV) is a near-space vehicle with broad prospects due to its high speed, large-scale maneuvering, and long-range gliding capabilities \cite{1}. In order to counter the missile defense system, it is necessary to enhance the penetration capability of the HGV. A practical penetration approach of the HGV is to make the best use of the speed superiority \cite{2}. However, the HGV may encounter multiple interceptors during the flight, and it is also essential to ensure that the HGV can successfully strike the target after a successful penetration. With the development of radar technology and anti-missile technology, especially cooperative interception technology \cite{3,4}, interceptors weaken the survival capability of the HGV. Therefore, the investigation of penetration trajectory design is essential for the HGV to confront interceptors.

The attack-defense confrontation in aerial engagement is essentially an optimal pursuit-evasion game. In \cite{7}, the optimal guidance and avoidance problem of missiles was solved by bounded control linear differential game, and the singular perturbation method was used to analyze the problem of medium-range air-to-air interception. Hagedorn et al. \cite{8} investigated a pursuit-evasion problem in which a faster evader must pass between two pursuers. The differential game of one evader and two pursuers with a nonconvex payoff function was investigated in \cite{9} and \cite{10}. However, for the HGV, an available penetration method is trajectory optimization. In \cite{11}, considering the coupled aerodynamics and thrust properties of hypersonic vehicles, a penetration trajectory programming method satisfying the saturation limit of control inputs was proposed using the improved model predictive static programming method. A cooperative guidance strategy for multiple HGVs systems was investigated in \cite{5}. A lateral pendulum maneuvering strategy was presented to improve the penetration performance of the HGV in \cite{12}. Based on the prerequisite penetration condition, a penetration trajectory optimization algorithm for an air-breathing hypersonic vehicle was designed in \cite{2}, which took the control costs as the objective function to minimize the fuel consumption and maneuver range. The penetration trajectory optimization of hypersonic vehicles considering the waypoint and no-fly zone constraints was discussed in \cite{6,13,14}. The penetration trajectory optimization of hypersonic vehicles striking a stationary target by a spiral-diving maneuver was studied in \cite{15}. Although there are many studies on HGV penetration, these studies mainly consider static or single-interceptor threats.  There are few existing results for the multi-interceptor penetration trajectory optimization, which is a practical and meaningful topic. Therefore, the HGV penetration method for more than one interceptor needs to be further focused on and investigated.

Since HGV is a highly nonlinear system, it is not always possible to obtain analytical solutions for nonlinear systems \cite{ALTAN2020,ALTAN2020con}. Numerical methods are suitable for dealing with highly nonlinear systems. It is also worth noting that real-time performance is essential in penetration trajectory optimization. Convex optimization, as a direct numerical method for trajectory optimization, is potential for online optimization \cite{LIU2021TAES}. Due to the excellent convergence performance and the broad applicability, the convex optimization method has been widely utilized in trajectory optimization \cite{ZHAO2021}. The trajectory optimization problem with complex constraints and nonlinear dynamics can be converted into a sequence of convex optimization problems by convexification techniques \cite{16,17,18,19,WANG2020linesearch,WANG2020linesearch2,20,21}. Convex optimization is widely used in complex trajectory optimization problems, such as Mars landing \cite{16}, maximum-crossrange problems \cite{19}, and atmospheric entry \cite{WANG2020linesearch2,20}. The convex optimization method can quickly converge to accurate solutions, and its potential application for real-time computational guidance has been wildly verified \cite{22,23,24,25,26,27,28,PEI2021realtime}. Convex optimization methods for solving non-convex optimal control problems in real-time were studied in \cite{23} using successive convexification approaches. In \cite{25}, real-time on-board trajectory generation for high-performance quad-rotor flight was investigated based on convexification techniques. Szmuk et al. \cite{26} proposed a real-time implementable successive convexification algorithm for the six-DOF powered descent guidance problem.  In \cite{27}, an onboard-implementable real-time algorithm based on convex optimization was presented for powered-descent guidance. A novel robust trajectory optimization procedure that combined the polynomial chaos theory with the convex optimization technique was proposed in \cite{29}. A minimum-fuel powered-descent optimal guidance algorithm that incorporates obstacle avoidance based on convex optimization was presented in \cite{30}. In \cite{31}, an advanced optimal-control optimizer based on successive convex optimization was developed and analyzed with a golden-section line-search method to enhance its convergence. A novel trajectory planning strategy for spacecraft relative motion was presented in \cite{32} by combining convex optimization and multi-resolution technique, and it performs better in computing efficiency than the traditional uniform discretization. Compared with the nonlinear programming-based method, the successive convex programming-based method might achieve faster convergence for trajectory optimization \cite{22}. Besides, it is also worth mentioning that machine-learning-based methods have a fast trajectory optimization and control ability \cite{ALTAN2018con}. For example, machine-learning methods have been used for entry guidance \cite{CHENG2021TAES} and spacecraft trajectory planning \cite{CHAI2020TNNLS}. However, the machine-learning method faces weak convergence property and hard reward function design \cite{LIU2021TAES}. Instead of relying on the quality of training data, the convex optimization methods can react to the environment because it actually calculates new trajectories \cite{Hofmann2021}. Thus, the convex optimization method is mainly investigated in the current study.

Motivated by the above research status and research difficulties, the penetration trajectory optimization of the HGV under the threat of two interceptors is studied. A mathematical model of the HGV penetration problem with multiple interceptors is established, and a penetration strategy is presented based on the initial line-of-sight (LOS) angle information of interceptors. The original nonlinear optimal control problem is transformed into second-order cone programming (SOCP) problems by convexification techniques. Then, the penetration trajectory is obtained by iterative solution of the successive SOCP method. The convergent trajectory obtained by the proposed method allows the HGV to evade two interceptors and reach the target area successfully. A variable trust region method is presented to improve the performance of the successive SOCP method. Compared with the line search approach, the proposed method can obtain the solution with similar optimality with less time consumption. The numerical simulation results show the effectiveness of the proposed method.

Compared with the existing work, the main contributions of this article are as follows. Firstly, the mathematical penetration model of the HGV encountering two interceptors is established based on the proposed penetration coordinate system. Within this model, penetration trajectory optimization method for one HGV encountering two interceptors can be studied while also considering the terminal target area. Although the penetration trajectory optimization methods were studied in \cite{11}, \cite{12} and \cite{13,14,15}, the specific threats of multi-interceptors were not investigated. Secondly, a new penetration strategy for the HGV encountering two interceptors is proposed, while the strategy for encountering only one interceptor was studied in \cite{2} and \cite{11}. The proposed strategy can increase the miss distance of both interceptors while ensuring the arrival of the target area. In addition, the proposed strategy only requires the initial LOS angle information of the two interceptors, while the guidance strategy of the interceptor is also required in \cite{2} and \cite{11}. Thirdly, a successive SOCP algorithm with the variable trust region is proposed, balancing the trade-off between time consumption and optimality. The line search strategy used in \cite{15} and \cite{19,WANG2020linesearch,WANG2020linesearch2} can improve the speed and robustness of convergence. However, the comparison of the numerical simulations shows that the proposed method performs better than line search in time consumption.

The rest of this paper is organized as follows. In Section \ref{section2}, the penetration trajectory optimization problem of the HGV is formulated. In Section \ref{section3}, the penetration scenario of the HGV is analyzed, and the penetration strategy for the HGV encountering two interceptors is proposed. In Section \ref{section4}, the penetration trajectory optimization problem is reformulated into a SOCP problem using existing convexification techniques. A successive SCOP method with a variable trust region is proposed to solve the SOCP problems. In Section \ref{section5}, numerical simulations are presented to demonstrate the effectiveness and performance of the penetration trajectory optimization method proposed in this paper. Section \ref{section6} summarizes the whole paper.

\section{Penetration problem with two interceptors}\label{section2}
In this section, the mathematical model of HGV penetration trajectory optimization is established. Then, a detailed formulation of the trajectory optimization problem for the HGV encountering two interceptors is presented.
\subsection{Penetration problem modeling}
\begin{figure}[thpb]
\vspace{-1em}
  \centering
  \includegraphics[height=6.0cm]{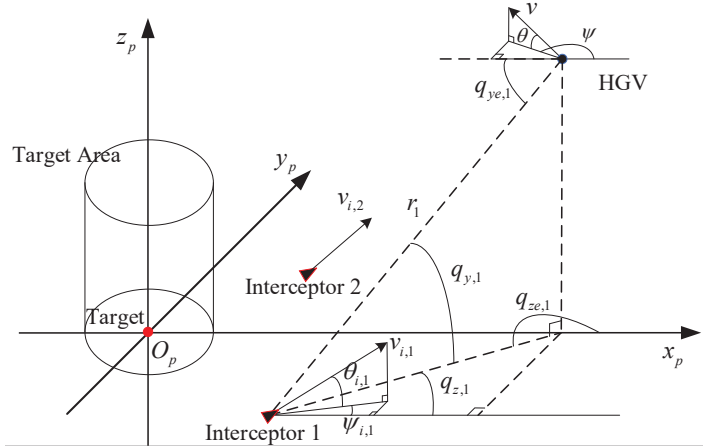}
  \vspace{-1em}\caption{Penetration geometry against two interceptors.}
  \label{fig1}
\end{figure}

The penetration trajectory optimization model presented in this subsection focuses on the gliding phase of the HGV. After gliding to the airspace near the stationary target, the HGV can track the diving trajectory or spiral-diving trajectory to strike the stationary target \cite{15}. In this paper, the trajectory optimization problem considering two interceptors and the terminal target area is investigated. Fig. \ref{fig1} shows the HGV penetration scenario in the gliding phase. The HGV needs to reach a target area after the penetration so that the vehicle can attack the stationary target effectively in the diving phase.

\begin{df}
\label{defination1}
The penetration coordinate system is a coordinate system with its origin ${O_p}$ at the stationary target. The axis ${O_p}{x_p}$ is located in the horizontal plane and points to the initial location of the HGV. The axis ${O_p}{z_p}$ is perpendicular to the horizontal plane and points upward. The axis ${O_p}{y_p}$ is perpendicular to the plane ${O_p}{x_p}{z_p}$ to complete a right-handed coordinate system.
\end{df}

Since the penetration process of the HGV encountering interceptors is mainly focused, the penetration trajectory optimization problem can be described as a short-range trajectory optimization problem. The curvature and rotation of the Earth are ignored in this paper. The three-dimensional motion equations of the HGV are given as follows
\begin{equation}\label{eq1}
\left\{ \begin{array}{l}
\dot h = v\sin \theta \\
{{\dot x}_p} = v\cos \theta \cos \psi \\
{{\dot y}_p} = v\cos \theta \sin \psi \\
\dot v =  - \frac{D}{m} - g\sin \theta \\
\dot \theta  = \frac{{L\cos \sigma }}{{mv}} - \frac{{g\cos \theta }}{v}\\
\dot \psi  = \frac{{L\sin \sigma }}{{mv\cos \theta }}
\end{array} \right.
\end{equation}
where $h$ is the altitude of the HGV, ${x_p}$ and ${y_p}$ represent the coordinate position of the vehicle on the ground, $v$ is the velocity of the vehicle, $\theta$ and $\psi$ are the flight path angle and the heading angle of velocity vector, $\alpha$ and $\sigma$ are the angle of attack and bank angle, $L$ and $D$ are the lift and drag, $g$ is the gravitational acceleration, $m$ is the vehicle mass. ${{q}_{ye}}$ and ${{q}_{ze}}$ are the elevation angle and the azimuth angle of the LOS from the HGV, respectively. Fig. \ref{fig1} shows explicit geometric meaning. ${{(\bullet )}_{,1}}$ and ${{(\bullet )}_{,2}}$ represent the corresponding variables of the two interceptors, respectively. ${r}_{i}$ is the relative distance between $i$th interceptor and the HGV , ${{q}_{y}}$ and ${{q}_{z}}$ are the elevation angle and the azimuth angle of the LOS from the interceptor, respectively. ${v}_{i}$ is the velocity of the interceptor, ${\theta}_{i}$ and ${\psi}_{i}$ are the flight path angle and the heading angle of the interceptor, respectively.

The aerodynamic force model of the HGV is given as follows \cite{2}
\begin{equation}\label{eq2}
\left\{ \begin{array}{l}
L = 0.5\rho {v^2}S({C_{L0}} + {C_{L\alpha }}\alpha )\\
D = 0.5\rho {v^2}S({C_{D0}} + {C_{D{\alpha ^2}}}{\alpha ^2})
\end{array} \right.
\end{equation}
where $\rho ={{\rho }_{0}}{{\operatorname{e}}^{\left( -h/{{h}_{s}} \right)}}$ is the atmospheric density, ${{\rho }_{0}}$ is the atmospheric density at sea level, ${{h}_{s}}=6700\text{m}$ is constant, $S$ is the aerodynamic reference area. The aerodynamic lift and drag coefficients are denoted as ${{C}_{L}}={{C}_{L0}}+{{C}_{L\alpha }}\alpha $ and ${{C}_{D}}={{C}_{D0}}+{{C}_{D{{\alpha }^{2}}}}{{\alpha }^{2}}$, respectively, where ${{C}_{L0}}$, ${{C}_{L\alpha }}$, ${{C}_{D0}}$, and ${{C}_{D{{\alpha }^{2}}}}$ are the aerodynamic parameters.

Due to the extensive range of bank angle and the direct numerical methods that may lead to the significant difference of bank angle in a single time step, the obtained trajectory may not conform to the physical characteristics of the HGV. Therefore, the bank angle rate is selected as a control variable \cite{22}. The range of the angle of attack is generally small, so the angle of attack is directly taken as another control variable. During the flight process, the angle of attack and bank angle rate are limited. The control constraints are given as follows
\begin{equation}\label{eq3}
\left\{ \begin{array}{l}
{\alpha _{\min }} \le \alpha  \le {\alpha _{\max }}\\
\left| {\dot \sigma } \right| \le {{\dot \sigma }_{\max }}
\end{array} \right.
\end{equation}
\begin{asmp}
\label{assumption1}
The initial velocity heading angle of the HGV is restrained as $\psi (t_0)\in \left( \pi /2,3\pi /2 \right)$.
\end{asmp}
\begin{rmk}
Assumption \ref{assumption1} is mild and reasonable. The HGV flies typically towards the target before encountering the interceptors. Thus, $\psi (t_0)\in \left( \pi /2,3\pi /2 \right)$ is reasonable in the penetration coordinate system.
\end{rmk}
\begin{asmp}
\label{assumption2}
During the flight of the HGV, the flight path angle and the heading angle are restrained as $\theta (t)\in \left( -\pi /2,\pi /2 \right)$ and $\psi (t)\in \left( \pi /2,3\pi /2 \right)$, respectively.
\end{asmp}
\begin{rmk}
Since $\theta (t) \in \left[ -\pi /2,\pi /2 \right]$ and the flight path angle of the HGV is generally impossible to reach $\pm \pi /2$. The flight path angle can meet Assumption \ref{assumption2} under normal circumstances. The HGV needs to reach the target area while penetrating the defense. That is, the HGV generally flies towards the target when encountering the interceptors. Thus, Assumption \ref{assumption2} is reasonable in the penetration coordinate system.
\end{rmk}
\begin{rmk}
The penetration coordinate system and the above assumptions play an essential role in the nonlinear dynamics transformation in Section \ref{section4}. In this paper, the penetration trajectory optimization problem can be solved by the SOCP solver after transformation in Section \ref{section4}.
\end{rmk}
\subsection{HGV penetration trajectory optimization problem}
The HGV needs to reach the target area. Since the target is the origin of the penetration coordinate system, the following objective function is defined
\begin{equation}\label{eq4}
{{J}_{1}}=\left\| ({{x}_{pf}},{{y}_{pf}}) \right\|
\end{equation}
where ${x}_{pf}$ and ${y}_{pf}$ are the terminal values of ${x}_{p}$ and ${y}_{p}$, respectively.

The HGV is expected to maneuver once the interceptors are detected. Then, the velocity direction should approach the expected penetration direction. Hence, the following objective function is defined
\begin{equation}\label{eq5}
{{J}_{2}}=\int\limits_{{{t}_{0}}}^{{{t}_{I}}}{\left\| (\theta ,\psi )-({{\theta }_{ex}},{{\psi }_{ex}}) \right\|dt}
\end{equation}
where ${{\theta }_{ex}}$ and ${{\psi }_{ex}}$ are the expected flight path angle and the expected heading angle, respectively. ${{\theta }_{ex}}$ and ${{\psi }_{ex}}$ will be presented in Section \ref{section3} based on the analysis of the penetration scenario.

The above optimal objective function (\ref{eq5}) indicates the expectation that the velocity direction of the HGV will approach the desired direction within the time $[{{t}_{0}},{{t}_{I}}]$. Since the vehicle has to reach the target area after penetration, (\ref{eq4}) and (\ref{eq5}) are combined into an objective function by adding a penalty term
\begin{equation}\label{eq6}
J = \left\| {({x_{pf}},{y_{pf}})} \right\| + {c_\vartheta }\int\limits_{{t_0}}^{{t_I}} {\left\| {(\theta ,\psi ) - ({\theta _{ex}},{\psi _{ex}})} \right\|dt}
\end{equation}
with constant ${{c}_{\vartheta }}>0$. In this way, the infeasibility problem is avoided \cite{17}.

Let $\boldsymbol{u}\text{= }\!\![\!\!\text{ }\alpha \text{,}\dot{\sigma }{{\text{ }\!\!]\!\!\text{ }}^{T}}$ denote the control, and $\boldsymbol{x}={{[h,{{x}_{p}},{{y}_{p}},v,\theta ,\psi ,\sigma ]}^{T}}$ denote the state. The state equations (\ref{eq1}) can be denoted as
\begin{equation}\label{eq7}
\boldsymbol{\dot{x}}=\boldsymbol{f}(\boldsymbol{x},\boldsymbol{u},t)
\end{equation}

The penetration trajectory optimization problem for the HGV is formulated as follows
\begin{equation}\label{eq8}
\begin{array}{l}
{\rm{P0}}:\min J = \left\| {({x_{pf}},{y_{pf}})} \right\| + {c_\vartheta }\int\limits_{{t_0}}^{{t_I}} {\left\| {(\theta ,\psi ) - ({\theta _{ex}},{\psi _{ex}})} \right\|dt} \\
{\rm{s}}{\rm{.t}}{\rm{.}}\left\{ \begin{array}{l}
{\boldsymbol{\dot x}} = {\boldsymbol{f}}({\boldsymbol{x}},{\boldsymbol{u}},t),\quad {\boldsymbol{x}}({t_0}) = {{\boldsymbol{x}}_0}\\
{\alpha _{\min }} \le \alpha (t) \le {\alpha _{\max }}\\
|\dot \sigma (t)| \le {{\dot \sigma }_{{\rm{max }}}}
\end{array} \right.
\end{array}
\end{equation}
where the purpose is to find the optimal control profile ${{\boldsymbol{u}}^{*}}={{[{{\alpha }^{*}},{{\dot{\sigma }}^{*}}]}^{T}}$ that minimizes the objective function (\ref{eq6}), and satisfies the dynamics (\ref{eq1}) and control constraints (\ref{eq3}).

\section{Penetration strategy}\label{section3}
In this section, the penetration scenario is analyzed. Then, the expected flight path angle and the expected heading angle are presented. The penetration strategy for one HGV encountering two interceptors is given.

\begin{figure}[htbp]
\centering
\subfigure[Head-on strike situation]{
\label {fig2-1}
\begin{minipage}[t]{0.33\linewidth}
\centering
\includegraphics[width=2in]{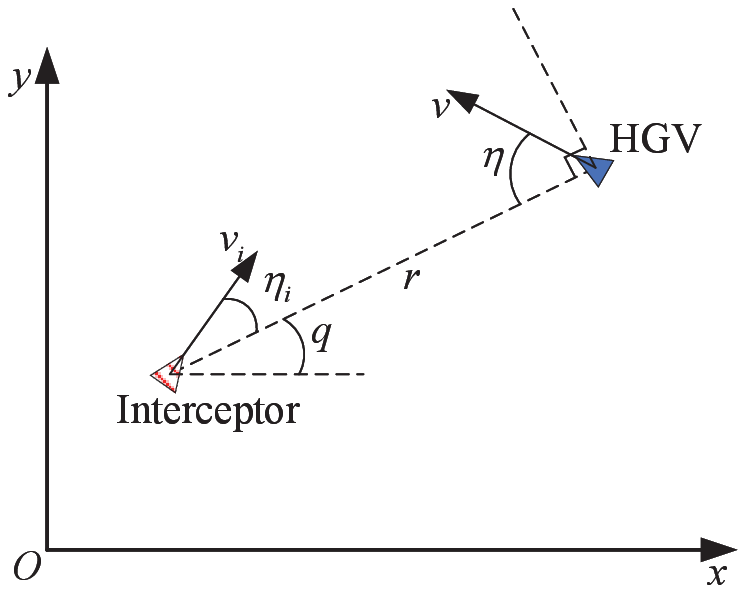}
\end{minipage}%
}%
\subfigure[Circumvent situation]{
\label {fig2-2}
\begin{minipage}[t]{0.33\linewidth}
\centering
\includegraphics[width=2in]{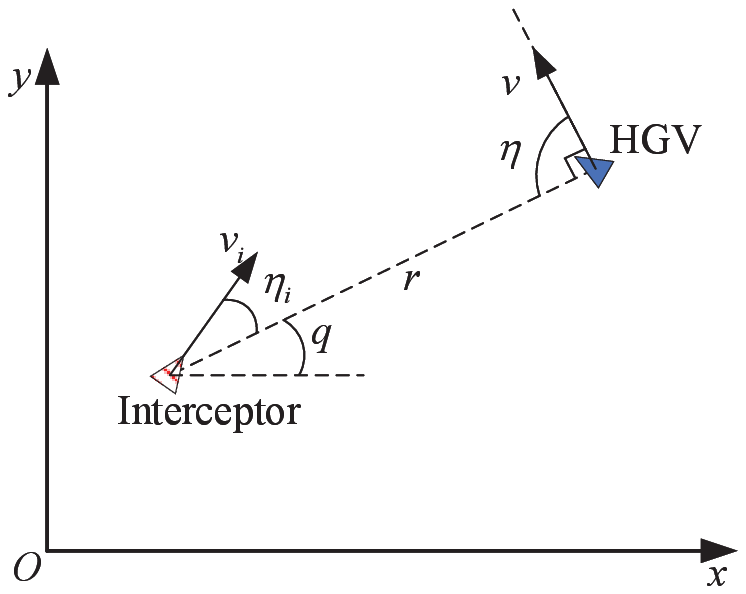}
\end{minipage}%
}%
\subfigure[Evasion situation]{
\label {fig2-3}
\begin{minipage}[t]{0.33\linewidth}
\centering
\includegraphics[width=2in]{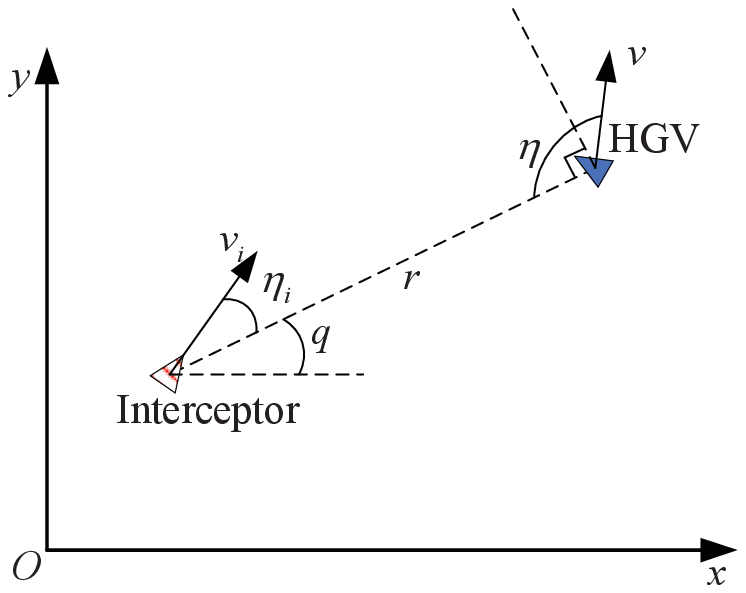}
\end{minipage}
}%
\centering
\caption{The two-dimensional penetration scenario.}
\label {fig2}
\end{figure}

The analysis starts from the two-dimensional penetration scenario confronting a single interceptor. Fig. \ref{fig2} shows three different situations in the two-dimensional penetration scenario.

The velocity of the HGV is usually much higher than that of the interceptors. Therefore, interceptors need to strike head-on to intercept the HGV effectively. The situations shown in Figs. \ref{fig2-2} and \ref{fig2-3} can be regarded as successful penetration. The geometry kinematics in the engagement can be given by
\begin{equation}\label{eq9}
\left\{ \begin{array}{l}
\dot r =  - v\cos \eta  - {v_i}\cos {\eta _i}\\
\dot qr = v\sin \eta  - {v_i}\sin {\eta _i}
\end{array} \right.
\end{equation}
where $r$ is the relative distance, $q$ is the LOS angle. $\eta $ is the angle between the velocity vector of the HGV and the LOS, and ${\eta}_{i} $ is the angle between the velocity vector of the interceptor and the LOS. The velocities of the HGV and the interceptor are considered constants in the penetration process analysis. As $\eta $ increases, $\dot{r}$ increases, which is beneficial to the HGV evasion. As $\eta $ increases, $\dot{q}$ also increases, which is beneficial for the HGV to circumvent the interceptor and penetrate. Therefore, with the increase of $\eta $, the penetration effectiveness of the HGV can be improved.

The purpose of the HGV is not just to evade but to circumvent the interceptor and reach the target area. When the velocity vector component in the direction perpendicular to the LOS angle is larger, the better the effectiveness of the HGV circumventing the interceptor will be. The velocity of the HGV is assumed to be constant in the analysis. Thus, the expected value of $\eta $ is $\pi /2$.

\begin{figure}[thpb]
\vspace{-1em}
  \centering
  \includegraphics[height=3.0cm]{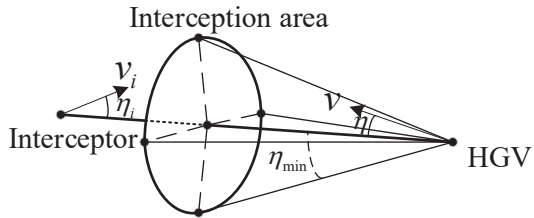}
  \vspace{-1em}\caption{The three-dimensional penetration scenario.}
  \label{fig3}
\end{figure}

According to the two-dimensional penetration scenario analysis, a similar analysis is carried out for the three-dimensional penetration scenario. Fig. \ref{fig3} shows the three-dimensional penetration scenario.
The geometric relationship of the penetration scenario can be obtained as
\begin{equation}\label{eq10}
\cos \eta =\cos {{q}_{ye}}\cos \theta \cos ({{q}_{ze}}-\psi )+\sin \theta \sin {{q}_{ye}}
\end{equation}
where ${{q}_{ye}}$ and ${{q}_{ze}}$ are the elevation angle and the azimuth angle of the LOS from the HGV, respectively.

If the velocity vector of the HGV is always outside the interception area of the interceptor, namely $\eta >{{\eta }_{\min }}$, the HGV can break through the head-on strike of the interceptor and achieve successful penetration \cite{2}. Thus, as angle $\eta $ increases, the penetration effectiveness will be improved.
\begin{df}
\label{definition2}
The penetration effectiveness of the HGV can be regarded as improved when $\eta $ increases.
\end{df}
\begin{thm}
\label{theorem1}
As $\left| {{q}_{ze}}-\psi  \right|$ and $\left| {{q}_{ye}}-\theta  \right|$ increase, the penetration effectiveness of the HGV can be improved.
\end{thm}
\proof
In order to analyze the HGV yaw motion for penetration, it is assumed that $\theta $ and ${{q}_{ye}}$ are constant. Since the cosine function is an even function, (\ref{eq10}) can be given as
\begin{equation}\label{eq11}
\cos \eta =\cos {{q}_{ye}}\cos \theta \cos (\left| {{q}_{ze}}-\psi  \right|)+\sin \theta \sin {{q}_{ye}}
\end{equation}
Take the partial derivative of (\ref{eq11}) with respect to $\left| {{q}_{ze}}-\psi  \right|$, it can be given as
\begin{equation}\label{eq12}
\frac{\partial (\cos \eta )}{\partial \left| {{q}_{ze}}-\psi  \right|}=-\cos {{q}_{ye}}\cos \theta \sin (\left| {{q}_{ze}}-\psi  \right|)
\end{equation}
Since the interceptor needs to intercept the HGV by head-on strike, ${{q}_{ye}}\in \left( -\pi /2,\pi /2 \right)$ and ${{q}_{ze}}\in \left( \pi /2,3\pi /2 \right)$ can be obtained. According to Assumption \ref{assumption2}, $\theta \in \left( -\pi /2,\pi /2 \right)$ and $\psi \in \left( \pi /2,3\pi /2 \right)$ are ensured. Thus, the following inequality can be established according to (\ref{eq12}).
\begin{equation}\label{eq13}
\frac{\partial (\cos \eta )}{\partial \left| {{q}_{ze}}-\psi  \right|}=-\cos {{q}_{ye}}\cos \theta \sin (\left| {{q}_{ze}}-\psi  \right|)<0
\end{equation}
It can be concluded that as $\left| {{q}_{ze}}-\psi  \right|$ increases, $\cos \eta $ decreases. Since $\eta  \in [0,\pi ]$, $\eta $ increases.

In order to analyze the HGV pitch motion for penetration, it is assumed that ${{q}_{ze}}=\psi $. In this case, the HGV flies toward the interceptor in the yaw direction. It is the most threatening situation to the HGV in the yaw direction. Therefore, the penetration strategy in the pitch direction can be analyzed when the HGV is most threatened in the yaw direction. (\ref{eq10}) can be rewritten as
\begin{equation}\label{eq14}
\cos \eta =\cos {{q}_{ye}}\cos \theta +\sin \theta \sin {{q}_{ye}}=\cos ({{q}_{ye}}-\theta )=\cos (\left| {{q}_{ye}}-\theta  \right|)
\end{equation}
Since ${{q}_{ye}}\in \left( -\pi /2,\pi /2 \right)$ and $\theta \in \left( -\pi /2,\pi /2 \right)$, as $\left| {{q}_{ye}}-\theta  \right|$ increases, $\cos \eta$ decreases.

It can be obtained that $\eta $ increases with the increase of $\left| {{q}_{ze}}-\psi  \right|$ and $\left| {{q}_{ye}}-\theta  \right|$. According to Definition \ref{definition2}, the penetration effectiveness of the HGV can be improved. The proof for Theorem \ref{theorem1} is completed. $\Box$

Similar to the two-dimensional penetration scenario analysis, the purpose of the HGV is not just to evade but to circumvent the interceptor and reach the target area. Therefore, when the HGV velocity direction is perpendicular to the LOS, that is, $\left| {{q}_{ye}}-{{\theta }} \right|=\pi /2$ and $\left| {{q}_{ze}}-{{\psi }} \right|=\pi /2$, the better the effectiveness of the HGV circumventing the interceptor.

According to the above analysis, a penetration strategy based on the expected flight path angle and the expected heading angle is proposed in this paper. The expected flight path angle is given by
\begin{equation}\label{eq15}
{\theta _{expect}}({q_{ye0}}) = \left\{ \begin{array}{l}
{q_{ye0}} + \frac{\pi }{2},{q_{ye0}} < 0\\
{q_{ye0}} - \frac{\pi }{2},{q_{ye0}} > 0\\
{q_{ye0}} + \frac{\pi }{2} - \chi ,{q_{ye0}} = 0
\end{array} \right.
\end{equation}
where $\chi $ is a positive small constant, which is introduced to ensure that the expected flight path angle also satisfies Assumption \ref{assumption2}. Similarly, the expected heading angle can be given by
\begin{equation}\label{eq16}
{\psi _{expect}}({q_{ze0}}) = \left\{ \begin{array}{l}
{q_{ze0}} + \frac{\pi }{2},{q_{ze0}} < \pi \\
{q_{ze0}} - \frac{\pi }{2},{q_{ze0}} > \pi \\
{q_{ze0}} + \frac{\pi }{2} - \chi ,{q_{ze0}} = \pi 
\end{array} \right.
\end{equation}
where ${{q}_{ye0}}$ and ${{q}_{ze0}}$ are the initial LOS elevation angle and initial LOS azimuth angle from the HGV, respectively. ${{\theta }_{0}}$ and ${{\psi }_{0}}$ are the initial flight path angle and initial heading angle, respectively.
\begin{thm}
\label{theorem2}
When the HGV encounters an interceptor coming from the front, that is, $\left| {{q}_{ye0}}-{{\theta }_{0}} \right|<\pi /2$ and $\left| {{q}_{ze0}}-{{\psi }_{0}} \right|<\pi /2$, if the flight path angle $\theta $ and the heading angle $\psi $ converge to the expected values given by (\ref{eq15}) and (\ref{eq16}), the penetration effectiveness of the HGV can be improved.
\end{thm}
\proof
If $\theta $ and $\psi $ converge to the expected values, it can be obtained from (\ref{eq15}) and (\ref{eq16}) that $\left| {{q}_{ye}}-\theta  \right|$ and $\left| {{q}_{ze}}-\psi  \right|$ converge to $\pi /2$. Thus, $\left| {{q_{ye}} - \theta } \right| > \left| {{q_{ye0}} - {\theta _0}} \right|$ and $\left| {{q_{ze}} - \psi } \right| > \left| {{q_{ze0}} - {\psi _0}} \right|$ are obtained. According to Theorem \ref{theorem1}, the penetration effectiveness of the HGV is improved. The proof for Theorem \ref{theorem2} is completed. $\Box$

Then, the penetration scenario for the multiple interceptors is discussed. In practical application, from the perspective of the efficiency-cost ratio, two interceptors are generally used to intercept one HGV. Thus, the scenario of two interceptors is mainly studied. The expected flight path angle and expected heading angle for two interceptors are given as follows based on the penetration strategy for a single interceptor (\ref{eq15}) and (\ref{eq16}).
\begin{equation}\label{eq17}
{\theta _{ex}} = \left\{ \begin{array}{l}
\max ({\theta _{expect}}({q_{ye0,1}}),{\theta _{expect}}({q_{ye0,2}})),{q_{ye0,1}} < 0,{q_{ye0,2}} < 0\\
\min ({\theta _{expect}}({q_{ye0,1}}),{\theta _{expect}}({q_{ye0,2}})),{q_{ye0,1}} > 0,{q_{ye0,2}} > 0\\
\frac{{{\theta _{expect}}({q_{ye0,1}}) + {\theta _{expect}}({q_{ye0,2}})}}{2},{q_{ye0,1}}{q_{ye0,2}} < 0\\
{q_{ye0,1}} - {\mathop{\rm sgn}} ({q_{ye0,2}})(\frac{\pi }{2} - \chi ),{q_{ye0,1}} = 0,{q_{ye0,2}} \ne 0\\
{q_{ye0,2}} - {\mathop{\rm sgn}} ({q_{ye0,1}})(\frac{\pi }{2} - \chi ),{q_{ye0,1}} \ne 0,{q_{ye0,2}} = 0\\
{\theta _{expect}}({q_{ye0,1}}),{q_{ye0,1}} = 0,{q_{ye0,2}} = 0
\end{array} \right.
\end{equation}
\begin{equation}\label{eq18}
{\psi _{ex}} = \left\{ \begin{array}{l}
\max ({\psi _{expect}}({q_{ze0,1}}),{\psi _{expect}}({q_{ze0,2}})),{q_{ze0,1}} < \pi ,{q_{ze0,2}} < \pi \\
\min ({\psi _{expect}}({q_{ze0,1}}),{\psi _{expect}}({q_{ze0,2}})),{q_{ze0,1}} > \pi ,{q_{ze0,2}} > \pi \\
\frac{{{\psi _{expect}}({q_{ze0,1}}) + {\psi _{expect}}({q_{ze0,2}})}}{2},({q_{ze0,1}} - \pi )({q_{ze0,2}} - \pi ) < 0\\
{q_{ze0,1}} - {\mathop{\rm sgn}} ({q_{ze0,2}} - \pi )(\frac{\pi }{2} - \chi ),{q_{ze0,1}} = \pi ,{q_{ze0,2}} \ne \pi \\
{q_{ze0,2}} - {\mathop{\rm sgn}} ({q_{ze0,1}} - \pi )(\frac{\pi }{2} - \chi ),{q_{ze0,1}} \ne \pi ,{q_{ze0,2}} = \pi \\
{\psi _{expect}}({q_{ze0,1}}),{q_{ze0,1}} = \pi ,{q_{ze0,2}} = \pi 
\end{array} \right.
\end{equation}
${{(\bullet )}_{,1}}$ and ${{(\bullet )}_{,2}}$ represent the corresponding variables of the two interceptors.
\begin{rmk}
The penetration strategy proposed in this paper is based on the initial LOS angle information of interceptors. Therefore, it is not necessary to obtain the guidance law and other information of interceptors. The HGV can begin to penetrate once the interceptors are discovered. From the perspective of the penetration scenario, the HGV should maneuver as soon as possible after detecting the interceptors. Thus, the penetration strategy based on the initial LOS angles of interceptors can ensure certain effectiveness. The proposed strategy will be analyzed and verified in the numerical simulations of Section \ref{section5}.
\end{rmk}

\section{Solution of penetration trajectory optimization problem}\label{section4}
The trajectory optimization problem P0 is a highly nonlinear optimal control problem with nonlinear dynamics and nonlinear objective function. Problem P0 can be reformulated into SOCP problems with the existing convexification techniques. The convexification process to transform the problem P0 into SOCP problems is presented in this section. A successive SOCP method with a variable trust region is developed to find solutions for the penetration trajectory optimization. A sequence of SOCP problems is iteratively solved until convergence is achieved.
\subsection{Linearization and Discretization of nonlinear dynamics}
In the penetration scenario, the terminal flight time of the HGV is challenging to determine, and the time is usually not constrained. Thus, the dynamics (\ref{eq1}) needs to be transformed. The derivation of ${{x}_{p}}$ leads to
\begin{equation}\label{eq19}
{{\dot{x}}_{p}}=v\cos \theta \cos \psi 
\end{equation}
Since $v$ is always positive. Under the Assumption \ref{assumption1} and Assumption \ref{assumption2},  $\cos \theta >0$ and $\cos \psi <0$ are ensured. Therefore, ${{x}_{p}}$ decreases monotonically with time, so it can be regarded as an independent variable of the state equations.

The dynamics (\ref{eq1}) can be transformed into
\begin{equation}\label{eq20}
\boldsymbol{{x}'}=\frac{d\boldsymbol{x}}{d{{x}_{p}}}=\boldsymbol{f}(\boldsymbol{x},\boldsymbol{u},{{x}_{p}})\text{,   }\boldsymbol{x}({{x}_{p0}})={{\boldsymbol{x}}_{0}}\text{   }
\end{equation}
where $\boldsymbol{x}={{[h,{{y}_{p}},v,\theta ,\psi ,\sigma ]}^{T}}$ is the state vector and the control vector is still $\boldsymbol{u}\text{= }\!\![\!\!\text{ }\alpha \text{,}\dot{\sigma }{{\text{ }\!\!]\!\!\text{ }}^{T}}$. In this paper, $\dot{\bullet }$ and ${\bullet }'$ represent derivatives of time and ${{x}_{p}}$, respectively. The state equations with ${{x}_{p}}$ as its independent variable can be obtained as follows.
\begin{equation}\label{eq21}
\left\{ \begin{array}{l}
\frac{{dh}}{{d{x_p}}} = \frac{{\tan \theta }}{{\cos \psi }}\\
\frac{{d{y_p}}}{{d{x_p}}} = \tan \psi \\
\frac{{dv}}{{d{x_p}}} =  - \frac{D}{{mv\cos \theta \cos \psi }} - \frac{{g\tan \theta }}{{v\cos \psi }}\\
\frac{{d\theta }}{{d{x_p}}} = \frac{{L\cos \sigma }}{{m{v^2}\cos \theta \cos \psi }} - \frac{g}{{{v^2}\cos \psi }}\\
\frac{{d\psi }}{{d{x_p}}} = \frac{{L\sin \sigma }}{{m{v^2}{{\cos }^2}\theta \cos \psi }}\\
\frac{{d\sigma }}{{d{x_p}}} = \frac{{\dot \sigma (t)}}{{v\cos \theta \cos \psi }}
\end{array} \right.
\end{equation}

One major work in the convexification of penetration problem lies in how to convexify the nonlinear dynamics. A convex constraint is obtained by successive linearizing the nonlinear dynamics about a reference solution $({{\boldsymbol{x}}^{(k)}},{{\boldsymbol{u}}^{(k)}})$, where superscripts $k$ is the current iteration.
\begin{equation}\label{eq22}
{\boldsymbol{x'}} = {\boldsymbol{f}}({\boldsymbol{x}},{\boldsymbol{u}},{x_p}) \approx A\left( {{{\boldsymbol{x}}^{(k)}},{{\boldsymbol{u}}^{(k)}},{x_p}} \right){\boldsymbol{x}} + B\left( {{{\boldsymbol{x}}^{(k)}},{{\boldsymbol{u}}^{(k)}},{x_p}} \right){\boldsymbol{u}} + {\boldsymbol{c}}\left( {{{\boldsymbol{x}}^{(k)}},{{\boldsymbol{u}}^{(k)}},{x_p}} \right)
\end{equation}
The full expressions of $A=\partial \boldsymbol{f}/\partial \boldsymbol{x}$, $B=\partial \boldsymbol{f}/\partial \boldsymbol{u}$, and $\boldsymbol{c}=\boldsymbol{f}(\boldsymbol{x},\boldsymbol{u},t)-A\boldsymbol{x}-B\boldsymbol{u}$ are reported in \ref{appendixa}. To ensure the convergence of successive linear approximations, the search space must be limited to a so-called trust region. The following trust region is given for validity of the linearization process.
\begin{equation}\label{eq23}
\left| \boldsymbol{x}-{{\boldsymbol{x}}^{(k)}} \right|\le {{\boldsymbol{\delta }}_{{0}}}
\end{equation}
${{\boldsymbol{\delta }}_{{0}}}\in {{\mathbb{R}}^{6}}$ is a constant vector. Note that the control constraints in (\ref{eq3}) can be regarded as a trust region $\left| \boldsymbol{u}-{{\boldsymbol{u}}^{(k)}} \right|\le {{\boldsymbol{\delta }}_{{u}}}$.

The problem is still an infinite-dimensional optimal control problem after successive linearization. Problem P0 must be properly discretized to facilitate the implementation of the successive SOCP method. The penetration process is averagely discretized through $N+1$ discretized points, and the penetration process is discretized into $N$ equal intervals. Thus, the infinite-dimensional optimal control problem is transformed into a finite-dimensional optimal control problem, imposing constraints on each discretized point. The step size is $\Delta {{x}_{p}}=({{x}_{pf}}-{{x}_{p0}})/N$, and the discretized points are denoted by $\{{{x}_{p0}},{{x}_{p1}},{{x}_{p2}},\cdots ,{{x}_{pN-1}},{{x}_{pN}}\}$. The corresponding state and control are discretized into ${{\boldsymbol{x}}_{i}}=\boldsymbol{x}({{x}_{pi}})$ and ${{\boldsymbol{u}}_{i}}=\boldsymbol{u}({{x}_{pi}})$, $i=0,1,\cdots ,N$. Then, the dynamics can be integrated numerically by the trapezoidal rule as follows \cite{17}.
\begin{equation}\label{eq24}
{{\boldsymbol{x}}_i} = {{\boldsymbol{x}}_{i - 1}} + \frac{{\Delta {x_p}}}{2}\left[ {(A_{i - 1}^{(k)}{{\boldsymbol{x}}_{i - 1}} + B_{i - 1}^{(k)}{{\boldsymbol{u}}_{i - 1}} + {\boldsymbol{c}}_{i - 1}^{(k)})} \right.\left. { + (A_i^{(k)}{{\boldsymbol{x}}_i} + B_i^{(k)}{{\boldsymbol{u}}_i} + {\boldsymbol{c}}_i^{(k)})} \right]{\rm{,    }}i = 0,1, \cdots ,N
\end{equation}
where $A_{i}^{(k)}=A\left( {{\boldsymbol{x}}^{(k)}}({{x}_{pi}}),{{\boldsymbol{u}}^{(k)}}({{x}_{pi}}),{{x}_{pi}} \right)$, $B_{i}^{(k)}=B\left( {{\boldsymbol{x}}^{(k)}}({{x}_{pi}}),{{\boldsymbol{u}}^{(k)}}({{x}_{pi}}),{{x}_{pi}} \right)$, and $\boldsymbol{c}_{i}^{(k)}=\boldsymbol{c}\left( {{\boldsymbol{x}}^{(k)}}({{x}_{pi}}),{{\boldsymbol{u}}^{(k)}}({{x}_{pi}}),{{x}_{pi}} \right)$. Let ${{H}_{i}}=-I+\frac{\Delta {{x}_{p}}}{2}A_{i}^{(k)}$, ${{H}_{i-1}}=I+\frac{\Delta {{x}_{p}}}{2}A_{i-1}^{(k)}$, and ${{G}_{i}}=\frac{\Delta {{x}_{p}}}{2}B_{i}^{(k)}$, (\ref{eq24}) is further rewritten as
\begin{equation}\label{eq25}
{H_{i - 1}}{{\boldsymbol{x}}_{i - 1}} + {H_i}{{\boldsymbol{x}}_i} + {G_{i - 1}}{{\boldsymbol{u}}_{i - 1}} + {G_i}{{\boldsymbol{u}}_i} =  - \frac{{\Delta {x_p}}}{2}({\boldsymbol{c}}_{i - 1}^{(k)} + {\boldsymbol{c}}_i^{(k)})
\end{equation}
where $I$ is a unit matrix with appropriate dimensions. The optimization variable vector can be denoted as $\boldsymbol{z}={{\left[ \begin{array}{*{35}{l}}
   \boldsymbol{x}_{0}^{T} & \cdots  & \boldsymbol{x}_{N}^{T} & \boldsymbol{u}_{0}^{T} & \cdots  & \boldsymbol{u}_{N}^{T}  \\
\end{array} \right]}^{T}}$. After successive linearization and discretization, the dynamics are transformed into a linear algebraic system in $\boldsymbol{z}$ as
\begin{equation}\label{eq26}
M{\boldsymbol{z}} = F
\end{equation}
where
\[M = \left[ {\begin{array}{*{20}{c}}
I&0&0& \cdots &0&0&0&0&0& \cdots &0&0\\
{{H_0}}&{{H_1}}&0& \cdots &0&0&{{G_0}}&{{G_1}}&0& \cdots &0&0\\
 \vdots & \vdots & \vdots & \ddots & \vdots & \vdots & \vdots & \vdots & \vdots & \ddots & \vdots & \vdots \\
0&0&0& \cdots &{{H_{N - 1}}}&{{H_N}}&0&0&0& \cdots &{{G_{N - 1}}}&{{G_N}}
\end{array}} \right]\]
\[F =  - \frac{{\Delta {x_p}}}{2}\left[ {\begin{array}{*{20}{c}}
{ - \left( {\frac{2}{{\Delta {x_p}}}} \right){{\boldsymbol{x}}_0}}\\
{{\boldsymbol{c}}_0^k + {\boldsymbol{c}}_1^k}\\
 \vdots \\
{{\boldsymbol{c}}_{N - 1}^k + {\boldsymbol{c}}_N^k}
\end{array}} \right]\]
It should be noted than the above equations also include initial state constraints $\boldsymbol{x}({{x}_{p0}})={{\boldsymbol{x}}_{0}}$.

\subsection{Convexification of the objective function}
In the penetration coordinate system, ${{x}_{p0}}>0$ can be ensured. When encountering interceptors, the primary goal of the HGV is to penetrate defenses. Set ${{x}_{pf}}={{x}_{pN}}=0$, the objective function in problem P0 can be transformed into
\begin{equation}\label{eq27}
J = \left\| {{y_{pf}}} \right\| + {c_\vartheta }\int\limits_{{t_0}}^{{t_I}} {\left\| {(\theta ,\psi ) - ({\theta _{ex}},{\psi _{ex}})} \right\|dt}
\end{equation}

Since the objective function can also be discretized, the above nonlinear objective function can be easily transformed into a linear function form by introducing two slack variables $\varsigma $ and ${{\vartheta }_{j}}$:
\begin{equation}\label{eq28}
J = \varsigma  + {c_\vartheta }\sum\limits_{j = 1}^{{N_I}} {{\vartheta _j}}
\end{equation}
subject to the following inequality constraints
\begin{equation}\label{eq29}
\left\| {{y_{pf}}} \right\| \le \varsigma {\rm{,       }}\left\| {({\theta _j},{\psi _j}) - ({\theta _{ex}},{\psi _{ex}})} \right\| \le {\vartheta _j}{\rm{,       }}j = 1,2, \cdots ,{N_I}
\end{equation}
The constraints in (\ref{eq29}) are second-order cone constraints, which are convex.

\subsection{SOCP form of penetration trajectory optimization problem}
The original problem P0 has been reformulated into a convex optimization problem P1
\begin{equation}\label{eq30}
{\rm{P1}}:\min J = \varsigma  + {c_\vartheta }\sum\limits_{j = 1}^{{N_I}} {{\vartheta _j}}
\end{equation}
subject to
\begin{equation}\label{eq31}
M\boldsymbol{z}=F
\end{equation}
\begin{equation}\label{eq32}
{{\alpha }_{\min }}\le {{\alpha }_{i}}\le {{\alpha }_{\max }}
\end{equation}
\begin{equation}\label{eq33}
|{{\dot{\sigma }}_{i}}|\le {{\dot{\sigma }}_{\text{max }}}
\end{equation}
\begin{equation}\label{eq34}
\left\| {{y}_{pf}} \right\|\le \varsigma \text{,       }\left\| ({{\theta }_{i}},{{\psi }_{i}})-({{\theta }_{ex}},{{\psi }_{ex}}) \right\|\le {{\vartheta }_{i}}\text{,       }i=1,2,\cdots ,{{N}_{I}}
\end{equation}
\begin{equation}\label{eq35}
\left| {{\boldsymbol{x}}_{i}}-\boldsymbol{x}_{i}^{(k)} \right|\le {\boldsymbol{{\delta }}_{0}}
\end{equation}
where $i=1,\cdots ,N$. When $i=1,2,\cdots ,{{N}_{I}}$, the HGV adopts the penetration strategy proposed in this paper to perform penetration maneuvers. When a small ${{N}_{I}}$ is selected, the duration of the penetration maneuver may be short. Thus, the penetration effectiveness may be limited. On the contrary, when a large ${{N}_{I}}$ is selected, the penetration effectiveness may be improved. But the HGV may not be able to accurately reach the target area.

In this article, ${N_I} = N/4$ is selected, which means that the penetration strategy is adopted by the HGV in the first quarter of the flight.
Since problem P1 has a linear cost function, all constraints are either linear constraints or second-order cone constraints. Problem P1 is a SOCP problem. While the discretization is sufficiently accurate, the solution of the problem P1 is a close approximation of the solution to problem P0 \cite{17}.

\subsection{Variable trust region successive SOCP method}
After obtaining the initial guess trajectory, in the subsequent iterations of convex optimization, the radius of the trust region determines the convergence performance of the sequential convex optimization \cite{32}. The solution procedure may converge slowly, and the main reason is that the dynamics are hard to satisfy. One possible solution is to tighten the trust-region constraint \cite{19}. If the radius of the trust region is too large, it may lead to a significant deviation from the original problem and difficult to converge. If the trust region is too small, the iterations may be limited, resulting in insufficient convergence. Therefore, to obtain satisfactory optimization results and improve the convergence performance of the method, a successive SOCP method with a variable trust region is designed. First, a trust-region radius ${{\boldsymbol{\delta }}_{0}}$ is set, and then the trust-region radius is tightened gradually with iterations. The variable trust region method can be given as follows.
\begin{equation}\label{eq36}
\left| {{\boldsymbol{x}}_{i}}-\boldsymbol{x}_{i}^{(k)} \right|\le \boldsymbol{\delta}
\end{equation}
\begin{equation}\label{eq37}
\boldsymbol{\delta} =\text{sgm}(k){\boldsymbol{\delta }_{0}}
\end{equation}
\begin{equation}\label{eq38}
\text{sgm}(k)=\frac{1}{1+{{e}^{\frac{k}{{{l}_{1}}}-{{l}_{2}}}}}
\end{equation}
where ${{l}_{1}}$ and ${{l}_{2}}$ are positive constants, $k$ is number of iterations. The sigmoid function is used in (\ref{eq38}). The variable trust region method can be adjusted by choosing ${{l}_{1}}$ and ${{l}_{2}}$. When $k<{{l}_{1}}{{l}_{2}}$, the trust-region tightening rate increases with the number of iterations. When $k>{{l}_{1}}{{l}_{2}}$, the trust-region tightening rate decreases as the number of iterations increases. The smaller the positive constant ${{l}_{1}}$ is, the faster the sigmoid function $\text{sgm}(k)$ converges to zero.

A successive SOCP method with a variable trust region is proposed to solve the penetration trajectory optimization problem. The successive SOCP method is detailed as follows.
\begin{enumerate}[Step 1:]
\item According to the initial LOS angles between the HGV and two interceptors, the expected heading angle and the expected flight path angle are calculated by (\ref{eq17}) and (\ref{eq18}).  
\item Set $k=0$. Choose an initial trajectory profile $({{\boldsymbol{x}}^{(0)}},{{\boldsymbol{u}}^{(0)}})$ and determine the parameters ${{l}_{1}}$ and ${{l}_{2}}$.  
\item At the $k+1$ times of iteration, problem P1 is setup by using the results of last iteration $({{\boldsymbol{x}}^{(k)}},{{\boldsymbol{u}}^{(k)}})$. The solution $({{\boldsymbol{x}}^{(k+1)}},{{\boldsymbol{u}}^{(k+1)}})$ can be obtained by solving problem P1.
\item Check if the stopping criteria (\ref{eq39}) is satisfied.
\begin{equation}\label{eq39}
\mathop {\max }\limits_i \left| {{\boldsymbol{x}}_i^{k + 1} - {\boldsymbol{x}}_i^k} \right| \le {\boldsymbol{\varepsilon }}
\end{equation}
where ${{\boldsymbol{\varepsilon }}}\in {{\mathbb{R}}^{6}}$ is a user-defined constant vector. If the convergence condition is satisfied, the optimal solution of original problem P0 is found to be $({{\boldsymbol{x}}^{(k+1)}},{{\boldsymbol{u}}^{(k+1)}})$, otherwise, go to Step5.
\item Set $k=k+1$ and update the trust region according to variable trust region strategy (\ref{eq37}). Then go back to Step2.
\end{enumerate}

In Step 2 of the above algorithm, the constant controls can be applied to the initial state, then the initial trajectory $({{\boldsymbol{x}}^{(0)}},{{\boldsymbol{u}}^{(0)}})$ is obtained by numerical integration of the dynamics. In the next section, the numerical simulations will show that the variable trust region is critical to ensure convergence performance of the algorithm.

\section{Numerical simulations}\label{section5}
In this section, some numerical simulations are presented to demonstrate the effectiveness and performance of the penetration trajectory optimization method proposed in this paper.

The initial values and parameters of simulations are listed in Table \ref{table1}. The initial trajectory $({{\boldsymbol{x}}^{(0)}},{{\boldsymbol{u}}^{(0)}})$ is obtained by giving the constant control guesses ${{\boldsymbol{u}}^{(0)}}={{[{{2}^{\circ }},0]}^{T}}$. In the variable trust region method, the parameter of (\ref{eq38}) is selected as ${{l}_{1}}=2.5$ and ${{l}_{2}}=5$ respectively. The constant in (\ref{eq30}) is selected as ${{c}_{\vartheta }}=1\times {{10}^{-6}}$.

The initial trust-region radius and the stopping criteria are given as follows
\[{\boldsymbol{\delta} _0} = {[{\rm{5000,5000,1000,}}\frac{{{\rm{40}}\pi }}{{180}}{\rm{,}}\frac{{{\rm{40}}\pi }}{{180}}{\rm{,}}\frac{{{\rm{40}}\pi }}{{180}}]^T}\]
\[{\boldsymbol{\varepsilon }} = {{\rm{[300,500,50,}}\frac{{{\rm{0}}{\rm{.5}}\pi }}{{180}}{\rm{,}}\frac{{{\rm{0}}{\rm{.5}}\pi }}{{180}}{\rm{,}}\frac{{{\rm{2}}\pi }}{{180}}{\rm{]}}^T}\]
The software package MOSEK is adopted to solve the SOCP problems iteratively. MOSEK is called by the YALMIP environment toolkit in MATLAB to solve the SOCP problems. The simulations are performed on a desktop computer with an Intel Core i9-10900K 3.70 GHz processor and 32 GB of RAM.

\begin{figure}[thpb]
\vspace{-1em}
  \centering
  \includegraphics[height=6.0cm]{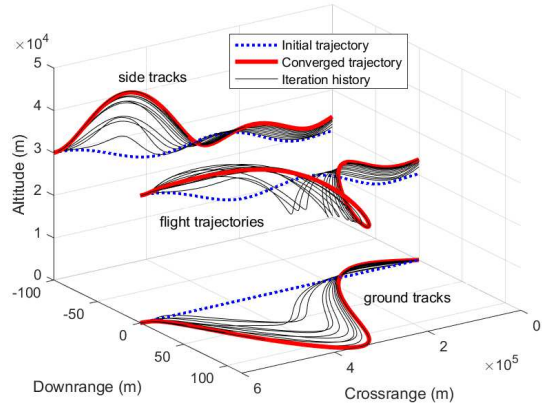}
  \vspace{-1em}\caption{Trajectories in each iteration.}
  \label{fig4}
\end{figure}

\begin{figure}[thpb]
\vspace{-1em}
  \centering
  \includegraphics[height=6.0cm]{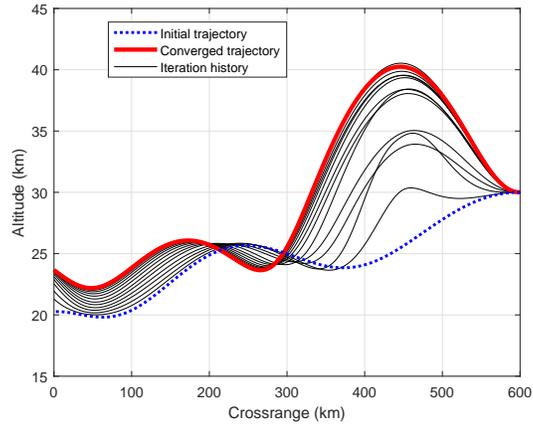}
  \vspace{-1em}\caption{Altitude profiles in each iteration.}
  \label{fig5}
\end{figure}

\begin{figure}[thpb]
\vspace{-1em}
  \centering
  \includegraphics[height=6.0cm]{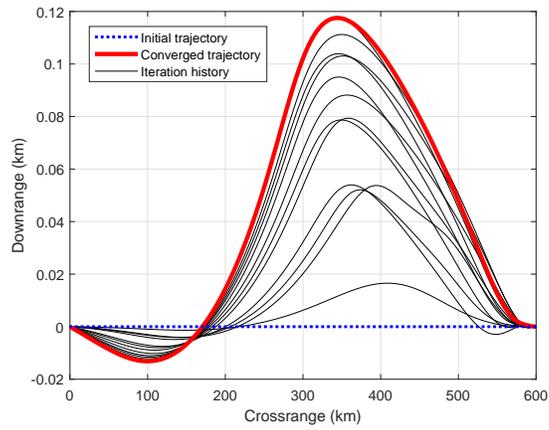}
  \vspace{-1em}\caption{Ground tracks in each iteration.}
  \label{fig6}
\end{figure}

\begin{figure}[thpb]
\vspace{-1em}
  \centering
  \includegraphics[height=6.0cm]{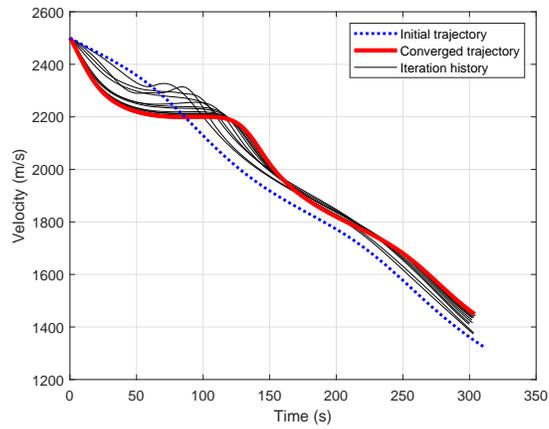}
  \vspace{-1em}\caption{Velocity profiles in each iteration.}
  \label{fig7}
\end{figure}

\begin{figure}[thpb]
\vspace{-1em}
  \centering
  \includegraphics[height=6.0cm]{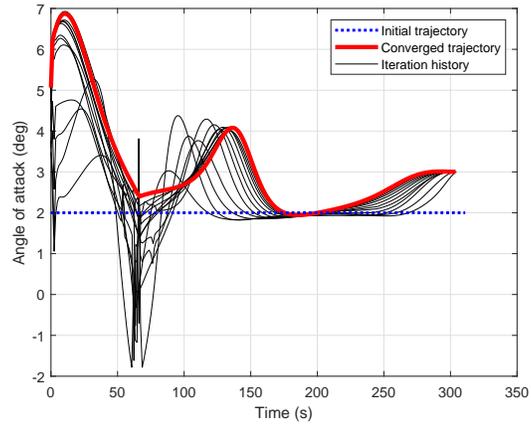}
  \vspace{-1em}\caption{Angle of attack profiles in each iteration.}
  \label{fig8}
\end{figure}

\begin{figure}[thpb]
\vspace{-1em}
  \centering
  \includegraphics[height=6.0cm]{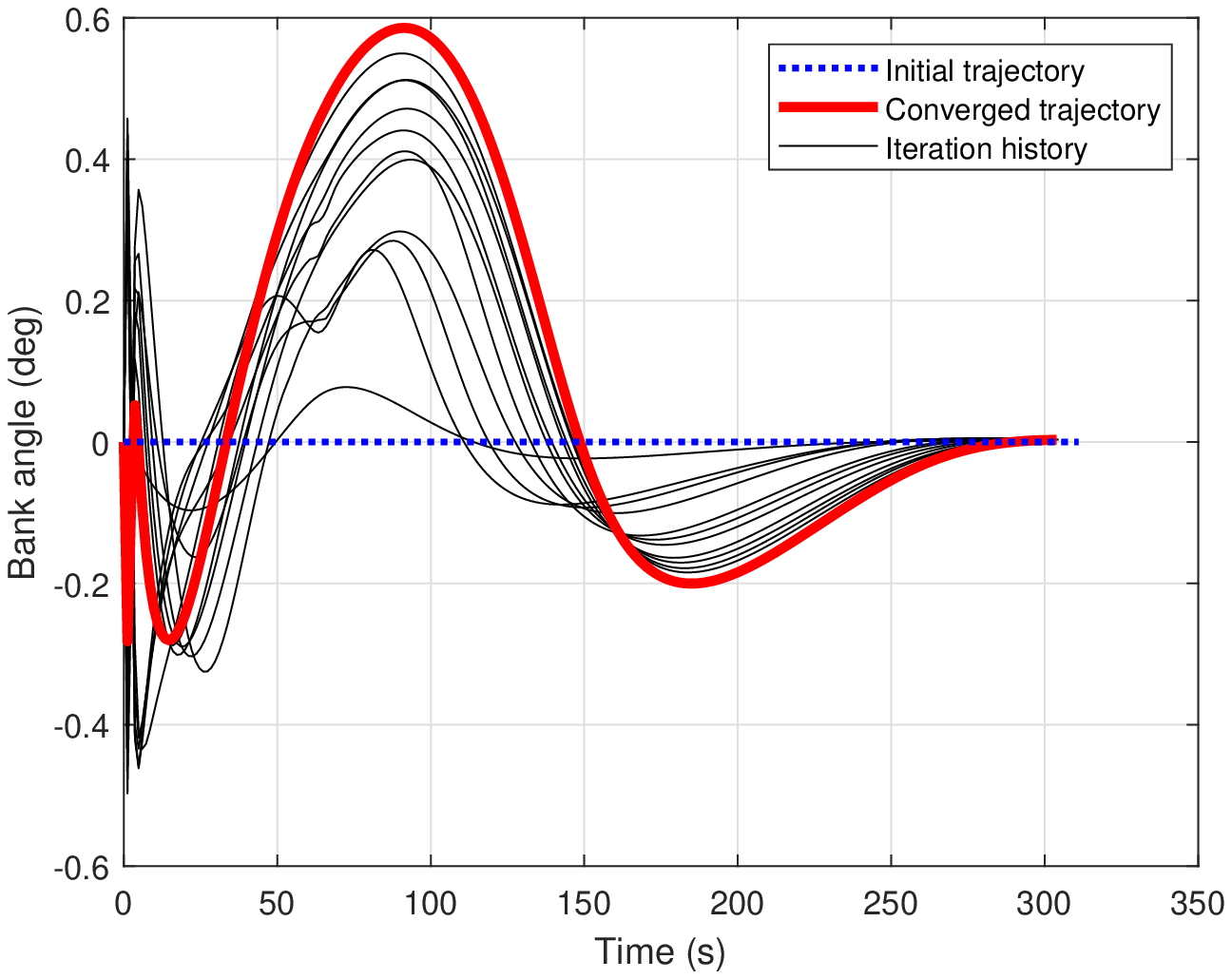}
  \vspace{-1em}\caption{Bank angle profiles in each iteration.}
  \label{fig9}
\end{figure}

\begin{figure}[thpb]
\vspace{-1em}
  \centering
  \includegraphics[height=6.0cm]{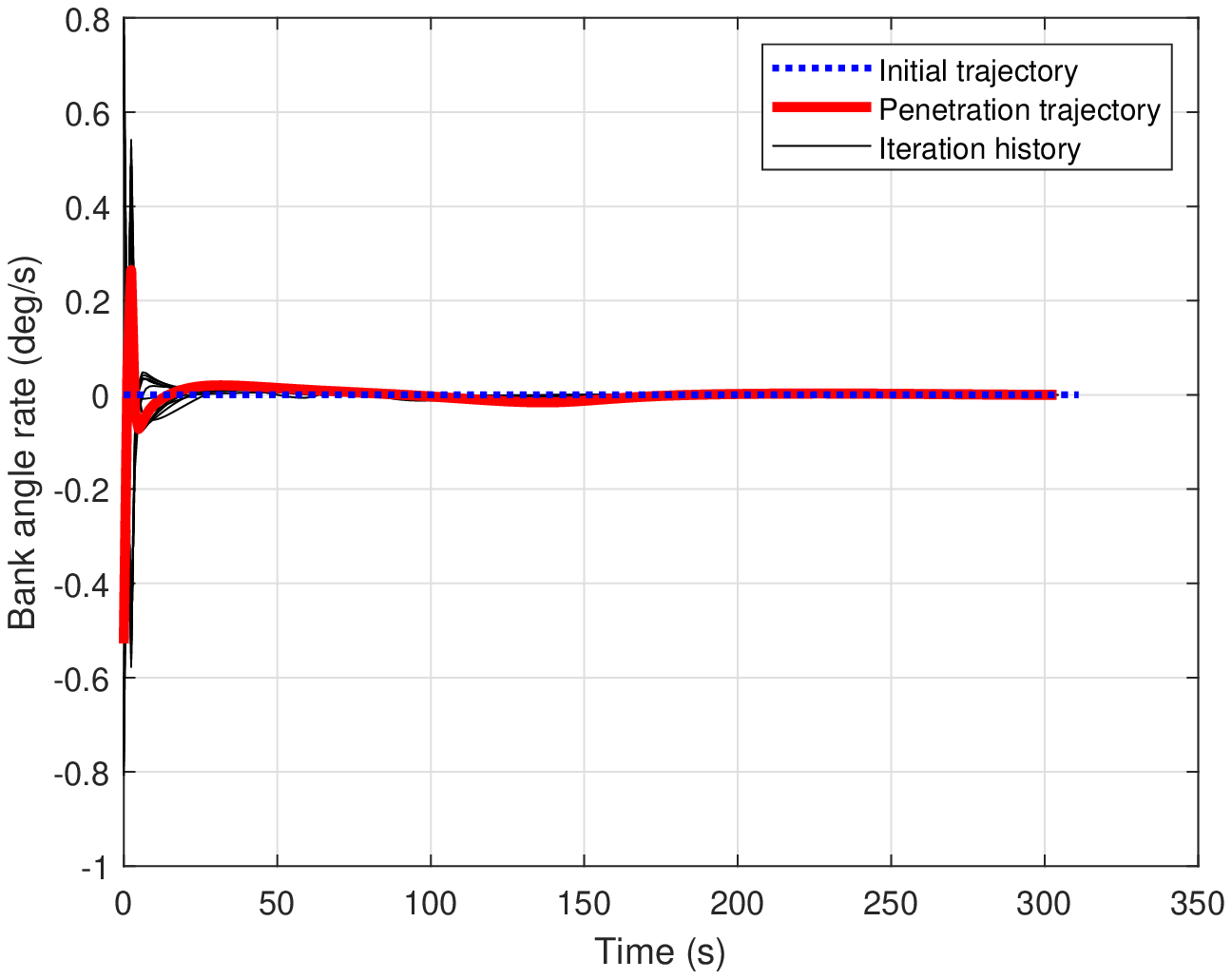}
  \vspace{-1em}\caption{Bank angle rate profiles in each iteration.}
  \label{fig10}
\end{figure}

\begin{figure}[thpb]
\vspace{-1em}
  \centering
  \includegraphics[height=6.0cm]{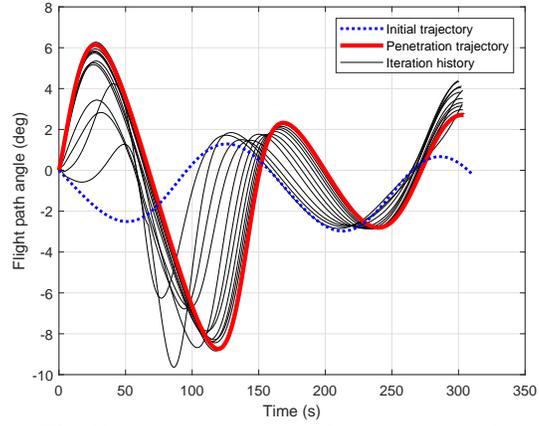}
  \vspace{-1em}\caption{Flight path angle profiles in each iteration.}
  \label{fig11}
\end{figure}

\begin{figure}[thpb]
\vspace{-1em}
  \centering
  \includegraphics[height=6.0cm]{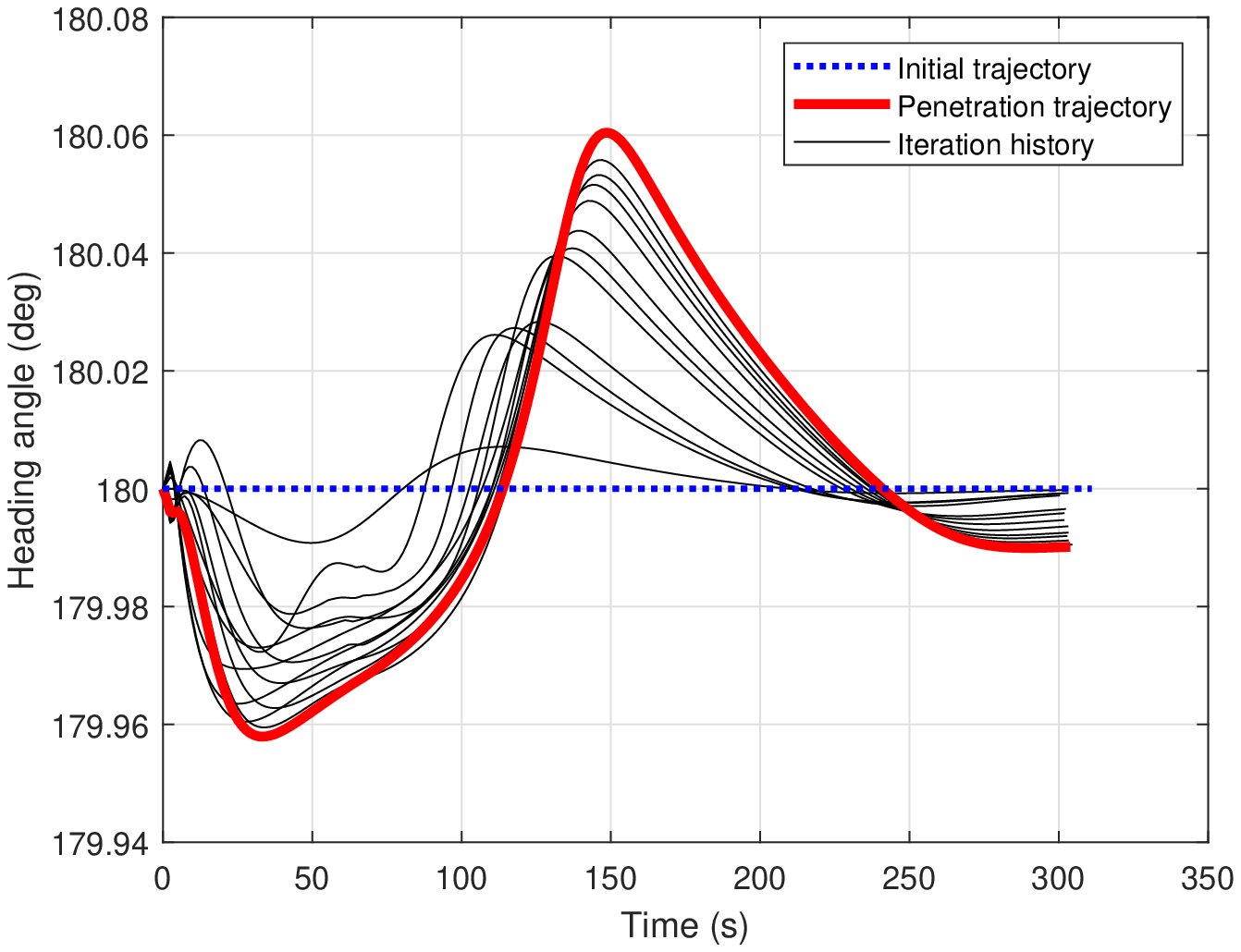}
  \vspace{-1em}\caption{Heading angle profiles in each iteration.}
  \label{fig12}
\end{figure}

By implementing the successive SOCP method with a variable trust region, a penetration trajectory is obtained in 14 iterations. In each iteration, MOSEK solves a SOCP problem in about 0.2s, and the total time of 14 iterations is 2.7188s. The results for all iterations are given in Figs. ~\ref{fig4}-\ref{fig12}. To further demonstrate the iteration process, the difference results of each state variable between consecutive iterations are reported in Table \ref{table2}. The difference of the states shown in Table \ref{table2} is defined as $|\Delta (\cdot )|:=\max \left| {{(\cdot )}^{(k)}}\left( {{t}_{i}} \right)-{{(\cdot )}^{(k-1)}}\left( {{t}_{i}} \right) \right|,i=1,2,\ldots ,N$, where $(\cdot )$ is each state variable.
\begin{table}
\centering
\caption{\newline Initial conditions of simulation.}
\label{table1}
\begin{tabular}{cc}
\hline
 Parameters &  Parameters Value\\
\hline
 $({{h}_{0}},{{x}_{p0}},{{y}_{p0}})$ (km) &  (30,600,0)\\
 ${{v}_{0}}$ (m/s) & 2500\\
 $({{\theta }_{0}},{{\psi }_{0}})$ (°) & (0,180)\\
 ${{\sigma }_{0}}$ (°) & 0\\
 $({{h}_{i0,1}},{{x}_{i0,1}},{{y}_{i0,1}})$ (km)  & (0,450, -5)\\
 $({{h}_{i0,2}},{{x}_{i0,2}},{{y}_{i0,2}})$ (km)  & (0,450,10)\\
 $N$  & 200\\
 ${{N}_{I}}$  & 50\\
 $m$ (kg)  & 802.2\\
 $({{C}_{L0}},{{C}_{L\alpha }})$  & (-0.013,1.833)\\
 $({{C}_{D0}},{{C}_{D{{\alpha }^{2}}}})$  & (0.015,4.596)\\
 $({{\alpha }_{\min }},{{\alpha }_{\max }})$ (°)  & (-4,10)\\
 ${{\dot{\sigma }}_{\text{max }}}$ (°/s)  & 5\\
\hline
\end{tabular}\\
\end{table}

\begin{table}
\centering
\caption{\newline Difference of the states between iterations.}
\label{table2}
\begin{tabular}{ccccccc}
\hline
 Iteration &  $\left| \Delta h \right|$ (km) & $\left| \Delta {{y}_{p}} \right|$ (m) & $\left| \Delta v \right|$ (m/s) & $\left| \Delta \theta  \right|$ (deg) & $\left| \Delta \psi  \right|$ (deg) & $\left| \Delta \sigma  \right|$ (deg/s) \\
\hline
1 & 4.5987  & 16.5573  & 84.1092  & 4.7894  & 0.0092 & 0.0967\\
2 & 4.6757 & 37.5858 & 54.7380 & 5.1188 & 0.0214 & 0.4733\\
3 & 1.8929 & 13.4552 & 43.1060 & 3.7305 & 0.0176 & 0.9548\\
4 & 1.5416 & 8.0203 & 19.6540 & 1.9641 & 0.0117 & 0.1806\\
5 & 3.5653 & 26.1087 & 45.5155 & 2.9722 & 0.0257 & 0.9011\\
6 & 1.2970 & 9.7387 & 20.9589 & 1.8689 & 0.0159 & 0.8519\\
7 & 0.7478 & 17.8217 & 10.3492 & 1.1056 & 0.0260 & 0.8129\\
8 & 1.0324 & 10.7013 & 15.8806 & 0.9645 & 0.0192 & 0.8392\\
9 & 0.5660 & 15.7127 & 7.8862 & 0.7980 & 0.0230 & 0.7909\\
10 & 0.3363 & 10.8558 & 8.1329 & 0.5554 & 0.0201 & 0.7648\\
11 & 0.3833 & 12.7165 & 5.9389 & 0.3865 & 0.0198 & 0.7383\\
12 & 0.6847 & 8.8389 & 12.0905 & 0.4046 & 0.0138 & 0.6963\\
13 & 0.3300 & 7.8281 & 8.4850 & 0.3484 & 0.0129 & 0.6111\\
14 & 0.0953 & 5.4188 & 1.3459 & 0.1657 & 0.0103 & 0.5219\\
\hline
\end{tabular}\\
\end{table}

The iterative process of penetration trajectory solution is shown in Figs. ~\ref{fig4}-\ref{fig7}. It can be seen clearly how the solution converges to the penetration trajectory from iteration to iteration. As shown in Fig. \ref{fig6}, the HGV can reach the target area. As shown in Figs. ~\ref{fig8}-\ref{fig10}, the controls can satisfy the control constraints. As shown in Figs. ~\ref{fig11}-\ref{fig12}, due to encountering interceptors from both sides, the HGV penetrates mainly through climbing maneuvers. Therefore, the bank angle and the bank angle rate are relatively small.

\begin{figure}[htbp]
\vspace{-1em}
  \centering
  \includegraphics[height=6.0cm]{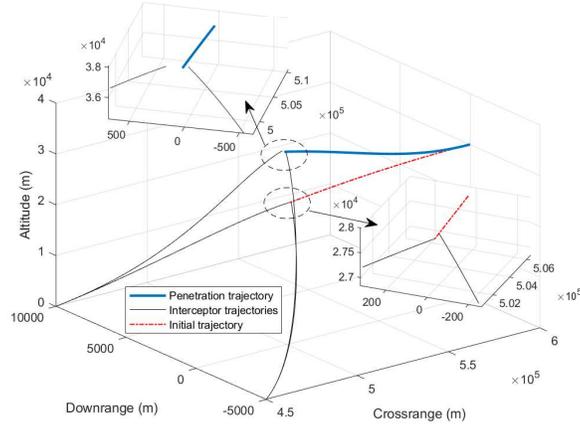}
  \vspace{-1em}\caption{Penetration detail.}
  \label{fig13}
\end{figure}

\begin{figure}[htbp]
\vspace{-1em}
  \centering
  \includegraphics[height=6.0cm]{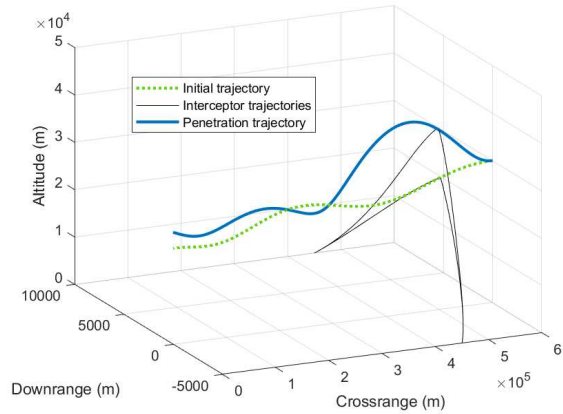}
  \vspace{-1em}\caption{Complete penetration process.}
  \label{fig14}
\end{figure}

To demonstrate the effectiveness of the penetration, a simulation is performed with two interceptors using the proportional navigation guidance. The guidance law constants of the interceptors are both selected as 5, and the maximum load of the interceptors is limited to 6g. The simulation step is set to 0.001s, according to the penetration trajectory, the final miss distances of interceptors are 41.939m and 150.8613m, respectively. When the HGV flies according to the initial trajectory, the final miss distances of the interceptors are 0.7622m and 1.6597m respectively. As shown in Figs. ~\ref{fig13}-\ref{fig14}, the HGV can achieve successful penetration by flying according to the penetration trajectory.


Furthermore, missions with various initial positions of interceptors from a given initial condition are also simulated to verify the effectiveness of the proposed method. The information for each mission is listed in Table \ref{table3}. The other conditions not listed in Table \ref{table3} are as same as Table \ref{table1}. The numbers of iterations indicate that the trajectory optimization problems of HGV penetration can be efficiently solved within a reasonable number of iterations. However, if the variable trust region is not used in the missions, the problem requires more than 40 iterations. While the number of iterations is quite different for each mission, the CPU time consumption of each iteration is almost the same. From the simulation results, the proposed method can successfully penetrate interceptors coming from different directions, which shows the effectiveness of the proposed penetration strategy. The last two columns in Table \ref{table3} imply that the final miss distances of mission 1 are less than that of mission 2 and 3. When encountering interceptors from the same side, the miss distances obtained by HGV are greater than those encountering interceptors from different sides.

\begin{figure}[htbp]
\vspace{-1em}
  \centering
  \includegraphics[height=6.0cm]{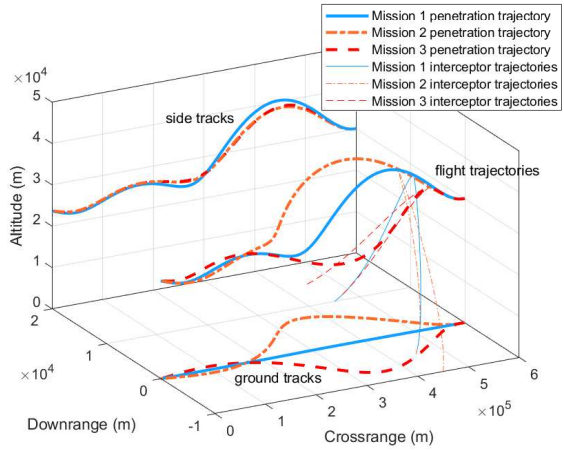}
  \vspace{-1em}\caption{Complete penetration process of missions 1-3.}
  \label{fig15}
\end{figure}

\begin{rmk}
As shown in Fig. \ref{fig15}, because the proposed penetration strategy only considers the initial LOS angle information of the interceptors, it also shows that the initial angle positions of the interceptors are significant to the interception efficiency. HGV mainly penetrates interceptors by longitudinal maneuver when encountering interceptors from different sides (i.e., mission 1). When encountering interceptors from the same direction (i.e., mission 2 and mission 3), HGV can effectively penetrate interceptors not only by longitudinal maneuver but also by lateral maneuver. The proposed method can be extended to the situation with more than two interceptors when HGV encounters interceptors with roughly the same attack direction. The HGV can penetrate through climbing maneuvers when more than two interceptors are launched from the ground. However, when too many interceptors come from different directions, the proposed penetration strategy may fail. The exploration of HGV encountering more than two interceptors, which prevent the HGV from penetrating any direction, is a possible direction of this research in the future.
\end{rmk}

\begin{table}
\centering
\caption{\newline Results of different missions.}
\label{table3}
\begin{threeparttable}
\begin{tabular}{ccccccc}
\hline
 Mission &  PI1\tnote{a} (km) & PI2\tnote{b} (km) & NI\tnote{c} & CTI\tnote{d} (s) & FMD1\tnote{e} (m) & FMD2\tnote{f} (m) \\
\hline
1 & (0,450,-5) & (0,450,10)  & 14  & 0.1942  & 41.9390 & 150.8613\\
2 & (0,450,-5) & (0,450,-10) & 13 & 0.1842 & 192.9452 & 351.9339\\
3 & (0,450,15) & (0,450,10) & 15 & 0.1673 & 492.1124 & 326.0877\\
\hline
\end{tabular}
      \begin{tablenotes}
        \footnotesize
        \item[a] Initial position of interceptor 1 $({{h}_{i0,1}},{{x}_{i0,1}},{{y}_{i0,1}})$.
	   \item[b] Initial position of interceptor 2 $({{h}_{i0,2}},{{x}_{i0,2}},{{y}_{i0,2}})$.
	\item[c] The number of iterations.
	\item[d] Average CPU time-consuming in each iteration.
	\item[e] Final miss distance of interceptor 1.
	\item[f] Final miss distance of interceptor 2.
      \end{tablenotes}
  \end{threeparttable}

\end{table}

\begin{figure}[htbp]
\vspace{-1em}
  \centering
  \includegraphics[height=6.0cm]{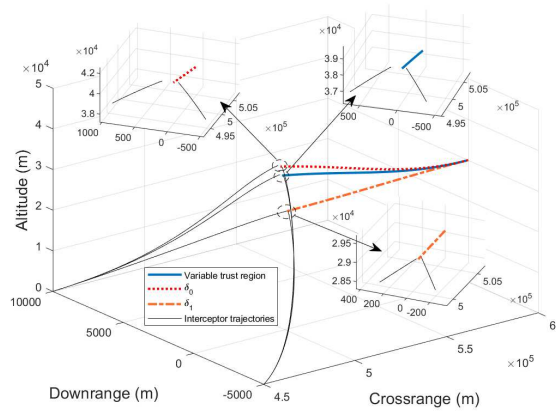}
  \vspace{-1em}\caption{Penetration detail of different trust region.}
  \label{fig16}
\end{figure}
\begin{figure}[htbp]
\vspace{-1em}
  \centering
  \includegraphics[height=6.0cm]{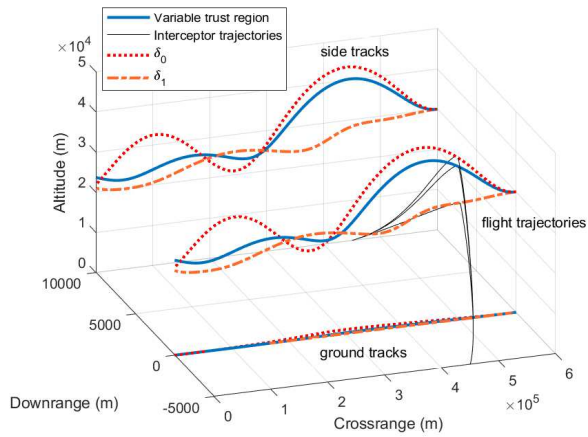}
  \vspace{-1em}\caption{Complete penetration process of different trust region.}
  \label{fig17}
\end{figure}

In order to demonstrate the effectiveness of the variable trust region method proposed in this paper, simulations are performed by giving two different constant trust regions. One of the constant trust regions is ${{\boldsymbol{\delta }}_{0}}$, which is a relatively large trust region. Another constant trust region is given below, which is a relatively small trust region.
\[{\boldsymbol{\delta }_1} = {[{\rm{2000,5000,500,}}\frac{{{\rm{20}}\pi }}{{180}}{\rm{,}}\frac{{{\rm{20}}\pi }}{{180}}{\rm{,}}\frac{{{\rm{20}}\pi }}{{180}}]^T}\]

The simulation results for different trust regions are shown in Figs. ~\ref{fig16}-\ref{fig17}. When the constant trust region ${{\boldsymbol{\delta }}_{1}}$ is used, the converged solution can be obtained by only 4 iterations. The total time of 4 iterations is only 0.7653s. However, the final miss distances of the interceptors are 2.9753m and 16.2671m, respectively. When a relatively small constant trust region is used, the feasible region is limited. The convergence is very fast, but the optimality is also limited. When the constant trust region ${{\boldsymbol{\delta }}_{0}}$ is used, the converged solution can be obtained by 87 iterations. The total time of 87 iterations is 16.6379s. The final miss distances of the interceptors are 59.9072m and 176.7311m, respectively. It is more time-consuming compared with the proposed method, and the miss distances are similar. Therefore, the method proposed in this paper can well balance the trade-off between time consumption and optimality.

\begin{rmk}
When the relatively small trust region ${{\boldsymbol{\delta }_1}}$ is selected, the feasible region of the optimization problem is limited to a small range. Therefore, the optimization problem converges quickly, but the convergent solution is limited to a small feasible region. However, when the relatively large trust region ${{\boldsymbol{\delta }}_{0}}$ is selected, the feasible region of the optimization problem is relatively large. It leads to the slow convergence of the optimization problem, but it can get a better solution with larger final miss distances. The proposed variable trust region method initially uses a relatively large feasible region ${{\boldsymbol{\delta }}_{0}}$, and the trust region gradually tightens with the iteration. The tightening rate accelerates gradually with the iteration. When the number of iterations reaches ${l}_{1}{l}_{2}$, the tightening rate slows down gradually. Table \ref{table2} and the related study \cite{22} show that the profile ${\boldsymbol{x}}^{(k)}$ changes significantly in the initial stage and will soon converge to a smaller range with the iteration.  Therefore, the proposed method, which gradually tightened the trust region, is reasonable. This study contributes to the understanding of how to make successive SOCP work for complex problems.
\end{rmk}

\begin{figure}[htbp]
\vspace{-1em}
  \centering
  \includegraphics[height=6.0cm]{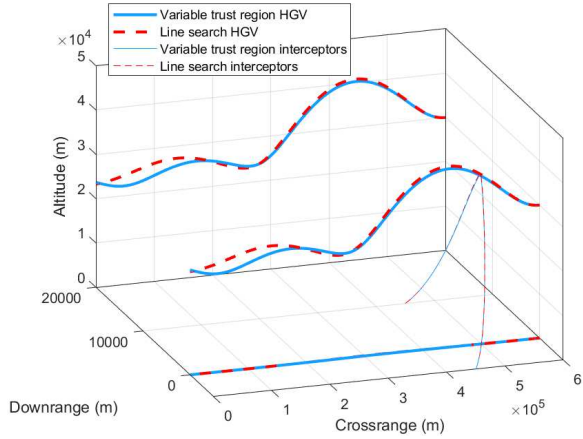}
  \vspace{-1em}\caption{Comparison of the trajectories.}
  \label{fig18}
\end{figure}

To further demonstrate the performance of the proposed variable trust region method, the line search approach \cite{19,WANG2020linesearch,WANG2020linesearch2}, a way to facilitate the convergence of the convex optimization, is used to solve the penetration trajectory optimization problem (i.e., mission 1) for comparison. The converged solution by our proposed method is presented and compared with the solution from the line search approach in Fig. \ref{fig18}. When using the line search approach, the converged solution can be obtained in 31 iterations.  The total time of 31 iterations is 5.6849s, and the final miss distances of the interceptors are 46.0676m and 160.6593m, respectively. The time consumption of the line search approach is nearly twice that of our method, but the final miss distances of the two methods are very close.

\begin{rmk}
When applying online, one possible way is to set a maximum iteration step \cite{PEI2021realtime}. Although it is not always possible to obtain the converged optimal solution, it can be considered that the obtained solution is near-optimal and will be updated in the following optimization circle. However, there is a trade-off between time consumption and optimality. Compared with the line search method, the simulation shows that the proposed method can obtain the solution with similar optimality with less time consumption. Therefore, our proposed method has great potential for real-time onboard applications.
\end{rmk}

\section{Conclusion}\label{section6}
In this article, the penetration trajectory optimization for HGV encountering two interceptors is proposed. The proposed method can optimize the trajectory so that the HGV can successfully penetrate the interceptors and reach the terminal target area. Based on the scenario analysis, a penetration strategy is proposed to improve penetration effectiveness, which only requires the initial LOS angle information of interceptors instead of the guidance law of interceptors and other information complicated to obtain in practice. Therefore, the proposed penetration strategy is conducive to engineering practice and application. The original nonconvex optimization problem is transformed into convex optimization problems by a series of technical methods. The transformed trajectory optimization problem can be solved efficiently by the existing convex optimization solver. The simulation results show that the proposed strategy can get a better penetration effectiveness encountering interceptors in the same attack direction, which shows that the initial angle positions of interceptors is significant for the interception. The results contribute to the design of not merely the penetration strategy but the interception strategy. 
Moreover, the proposed variable trust region method can well balance the trade-off between time consumption and optimality. Compared with the line search approach, our method consumes less time when the optimality is similar, showing the potential of the real-time application. Real-time trajectory optimization is a topic for further investigations. HGV against more than two interceptors attacking in all directions is also a possible direction of this research in the future.

\appendix
\section{Linearization Matrices}
\label{appendixa}
The analytical expressions of the non-zero elements of the $A$ matrix are detailed as follows.
\begin{equation}
{{A}_{14}}=\frac{1}{\cos \psi {{(\cos \theta )}^{2}}}
\end{equation}
\begin{equation}
{{A}_{15}}=\frac{\sin \psi \tan \theta }{{{(\cos \psi )}^{2}}}
\end{equation}
\begin{equation}
{{A}_{25}}=\frac{1}{{{(\cos \psi )}^{2}}}
\end{equation}
\begin{equation}
{{A}_{31}}=\frac{{{\rho }_{0}}{{e}^{-\frac{h}{{{h}_{s}}}}}vS({{C}_{D0}}+{{C}_{D{{\alpha }^{2}}}}{{\alpha }^{2}})}{2{{h}_{s}}m\cos \psi \cos \theta }
\end{equation}
\begin{equation}
{{A}_{33}}=\frac{\frac{{{\rho }_{0}}{{e}^{-\frac{h}{{{h}_{s}}}}}{{v}^{2}}S({{C}_{D0}}+{{C}_{D{{\alpha }^{2}}}}{{\alpha }^{2}})}{2m}+g\sin \theta }{{{v}^{2}}\cos \psi \cos \theta }-\frac{{{\rho }_{0}}{{e}^{-\frac{h}{{{h}_{s}}}}}S({{C}_{D0}}+{{C}_{D{{\alpha }^{2}}}}{{\alpha }^{2}})}{m\cos \psi \cos \theta }
\end{equation}
\begin{equation}
{{A}_{34}}=-\frac{g}{v\cos \psi }-\frac{\sin \theta (\frac{{{\rho }_{0}}{{e}^{-\frac{h}{{{h}_{s}}}}}{{v}^{2}}S({{C}_{D0}}+{{C}_{D{{\alpha }^{2}}}}{{\alpha }^{2}})}{2m}+g\sin \theta )}{v\cos \psi {{(\cos \theta )}^{2}}}
\end{equation}
\begin{equation}
{{A}_{35}}=-\frac{\sin \psi (\frac{{{\rho }_{0}}{{e}^{-\frac{h}{{{h}_{s}}}}}{{v}^{2}}S({{C}_{D0}}+{{C}_{D{{\alpha }^{2}}}}{{\alpha }^{2}})}{2m}+g\sin \theta )}{v{{(\cos \psi )}^{2}}\cos \theta }
\end{equation}
\begin{equation}
{{A}_{41}}=-\frac{{{\rho }_{0}}{{e}^{-\frac{h}{{{h}_{s}}}}}S({{C}_{L0}}+{{C}_{L\alpha }}\alpha )\cos \sigma }{2{{h}_{s}}m\cos \psi \cos \theta }
\end{equation}
\begin{equation}
{{A}_{43}}=\frac{\frac{g\cos \theta }{v}-\frac{{{\rho }_{0}}{{e}^{-\frac{h}{{{h}_{s}}}}}vS({{C}_{L0}}+{{C}_{L\alpha }}\alpha )\cos \sigma }{2m}}{{{v}^{2}}\cos \psi \cos \theta }+\frac{\frac{g\cos \theta }{{{v}^{2}}}-\frac{{{\rho }_{0}}{{e}^{-\frac{h}{{{h}_{s}}}}}S({{C}_{L0}}+{{C}_{L\alpha }}\alpha )\cos \sigma }{2m}}{v\cos \psi \cos \theta }
\end{equation}
\begin{equation}
{{A}_{44}}=\frac{g\sin \theta }{{{v}^{2}}\cos \psi \cos \theta }-\frac{\sin \theta (\frac{g\cos \theta }{v}-\frac{{{\rho }_{0}}{{e}^{-\frac{h}{{{h}_{s}}}}}vS({{C}_{L0}}+{{C}_{L\alpha }}\alpha )\cos \sigma }{2m})}{v\cos \psi {{(\cos \theta )}^{2}}}
\end{equation}
\begin{equation}
{{A}_{45}}=-\frac{\sin \psi (\frac{g\cos \theta }{v}-\frac{{{\rho }_{0}}{{e}^{-\frac{h}{{{h}_{s}}}}}vS({{C}_{L0}}+{{C}_{L\alpha }}\alpha )\cos \sigma }{2m})}{v{{(\cos \psi )}^{2}}\cos \theta }
\end{equation}
\begin{equation}
{{A}_{46}}=-\frac{{{\rho }_{0}}{{e}^{-\frac{h}{{{h}_{s}}}}}S({{C}_{L0}}+{{C}_{L\alpha }}\alpha )\sin \sigma }{2m\cos \psi \cos \theta }
\end{equation}
\begin{equation}
{{A}_{51}}=-\frac{{{\rho }_{0}}{{e}^{-\frac{h}{{{h}_{s}}}}}S({{C}_{L0}}+{{C}_{L\alpha }}\alpha )\sin \sigma }{2{{h}_{s}}m\cos \psi {{(\cos \theta )}^{2}}}
\end{equation}
\begin{equation}
{{A}_{54}}=\frac{{{\rho }_{0}}{{e}^{-\frac{h}{{{h}_{s}}}}}S({{C}_{L0}}+{{C}_{L\alpha }}\alpha )\sin \sigma \sin \theta }{m\cos \psi {{(\cos \theta )}^{3}}}
\end{equation}
\begin{equation}
{{A}_{55}}=\frac{{{\rho }_{0}}{{e}^{-\frac{h}{{{h}_{s}}}}}S({{C}_{L0}}+{{C}_{L\alpha }}\alpha )\sin \sigma \sin \psi }{2m{{(\cos \psi )}^{2}}{{(\cos \theta )}^{2}}}
\end{equation}
\begin{equation}
{{A}_{56}}=\frac{{{\rho }_{0}}{{e}^{-\frac{h}{{{h}_{s}}}}}S({{C}_{L0}}+{{C}_{L\alpha }}\alpha )\cos \sigma }{2m\cos \psi {{(\cos \theta )}^{2}}}
\end{equation}
\begin{equation}
{{A}_{63}}=-\frac{{\dot{\sigma }}}{{{v}^{2}}\cos \psi \cos \theta }
\end{equation}
\begin{equation}
{{A}_{64}}=-\frac{\dot{\sigma }\sin \theta }{v\cos \psi {{(\cos \theta )}^{2}}}
\end{equation}
\begin{equation}
{{A}_{65}}=-\frac{\dot{\sigma }\sin \psi }{v{{(\cos \psi )}^{2}}\cos \theta }
\end{equation}

The $B$ matrix is detailed as follows.
\begin{equation}
B=\left[ \begin{matrix}
   0 & 0  \\
   0 & 0  \\
   \frac{-{{\rho }_{0}}{{e}^{-\frac{h}{{{h}_{s}}}}}vS{{C}_{D{{\alpha }^{2}}}}\alpha }{m\cos \psi \cos \theta } & 0  \\
   \frac{{{\rho }_{0}}{{e}^{-\frac{h}{{{h}_{s}}}}}S{{C}_{L\alpha }}\cos \sigma }{2m\cos \psi \cos \theta } & 0  \\
   \frac{{{\rho }_{0}}{{e}^{-\frac{h}{{{h}_{s}}}}}S{{C}_{L\alpha }}\sin \sigma }{2m\cos \psi {{(\cos \theta )}^{2}}} & 0  \\
   0 & \frac{1}{v\cos \psi \cos \theta }  \\
\end{matrix} \right]
\end{equation}

The $\boldsymbol{c}$ vector elements are as follows.
\begin{equation}
{{c}_{1}}=\frac{\tan \theta }{\cos \psi }+\frac{\theta }{\cos \psi {{(\cos \theta )}^{2}}}+\frac{\psi \sin \psi \tan \theta }{{{(\cos \psi )}^{2}}}
\end{equation}
\begin{equation}
{{c}_{2}}=\tan \psi +\frac{\psi }{{{(\cos \psi )}^{2}}}
\end{equation}
\begin{equation}
\begin{array}{l}
{c_3} =  - \frac{D}{{mv\cos \theta \cos \psi }} - \frac{{g\tan \theta }}{{v\cos \psi }} + \frac{{h{\rho _0}{e^{ - \frac{h}{{{h_s}}}}}vS({C_{D0}} + {C_{D{\alpha ^2}}}{\alpha ^2})}}{{2{h_s}m\cos \psi \cos \theta }}\\
{\rm{       }} + v(\frac{{\frac{{{\rho _0}{e^{ - \frac{h}{{{h_s}}}}}{v^2}S({C_{D0}} + {C_{D{\alpha ^2}}}{\alpha ^2})}}{{2m}} + g\sin \theta }}{{{v^2}\cos \psi \cos \theta }} - \frac{{{\rho _0}{e^{ - \frac{h}{{{h_s}}}}}S({C_{D0}} + {C_{D{\alpha ^2}}}{\alpha ^2})}}{{m\cos \psi \cos \theta }})\\
{\rm{       }} + \theta ( - \frac{g}{{v\cos \psi }} - \frac{{\sin \theta (\frac{{{\rho _0}{e^{ - \frac{h}{{{h_s}}}}}{v^2}S({C_{D0}} + {C_{D{\alpha ^2}}}{\alpha ^2})}}{{2m}} + g\sin \theta )}}{{v\cos \psi {{(\cos \theta )}^2}}})\\
{\rm{       }} + \psi ( - \frac{{\sin \psi (\frac{{{\rho _0}{e^{ - \frac{h}{{{h_s}}}}}{v^2}S({C_{D0}} + {C_{D{\alpha ^2}}}{\alpha ^2})}}{{2m}} + g\sin \theta )}}{{v{{(\cos \psi )}^2}\cos \theta }})\\
{\rm{       }} + \alpha (\frac{{ - {\rho _0}{e^{ - \frac{h}{{{h_s}}}}}vS{C_{D{\alpha ^2}}}\alpha }}{{m\cos \psi \cos \theta }})
\end{array}
\end{equation}
\begin{equation}
\begin{array}{l}
{c_4} = \frac{{L\cos \sigma }}{{m{v^2}\cos \theta \cos \psi }} - \frac{g}{{{v^2}\cos \psi }} - \frac{{h{\rho _0}{e^{ - \frac{h}{{{h_s}}}}}S({C_{L0}} + {C_{L\alpha }}\alpha )\cos \sigma }}{{2{h_s}m\cos \psi \cos \theta }}\\
{\rm{       }} + v(\frac{{\frac{{g\cos \theta }}{v} - \frac{{{\rho _0}{e^{ - \frac{h}{{{h_s}}}}}vS({C_{L0}} + {C_{L\alpha }}\alpha )\cos \sigma }}{{2m}}}}{{{v^2}\cos \psi \cos \theta }} + \frac{{\frac{{g\cos \theta }}{{{v^2}}} - \frac{{{\rho _0}{e^{ - \frac{h}{{{h_s}}}}}S({C_{L0}} + {C_{L\alpha }}\alpha )\cos \sigma }}{{2m}}}}{{v\cos \psi \cos \theta }})\\
{\rm{       }} + \theta (\frac{{g\sin \theta }}{{{v^2}\cos \psi \cos \theta }} - \frac{{\sin \theta (\frac{{g\cos \theta }}{v} - \frac{{{\rho _0}{e^{ - \frac{h}{{{h_s}}}}}vS({C_{L0}} + {C_{L\alpha }}\alpha )\cos \sigma }}{{2m}})}}{{v\cos \psi {{(\cos \theta )}^2}}})\\
{\rm{       }} + \psi ( - \frac{{\sin \psi (\frac{{g\cos \theta }}{v} - \frac{{{\rho _0}{e^{ - \frac{h}{{{h_s}}}}}vS({C_{L0}} + {C_{L\alpha }}\alpha )\cos \sigma }}{{2m}})}}{{v{{(\cos \psi )}^2}\cos \theta }})\\
{\rm{       }} - \frac{{\sigma {\rho _0}{e^{ - \frac{h}{{{h_s}}}}}S({C_{L0}} + {C_{L\alpha }}\alpha )\sin \sigma }}{{2m\cos \psi \cos \theta }}\\
{\rm{       }} + \alpha (\frac{{{\rho _0}{e^{ - \frac{h}{{{h_s}}}}}S{C_{L\alpha }}\cos \sigma }}{{2m\cos \psi \cos \theta }})
\end{array}
\end{equation}
\begin{equation}
\begin{array}{l}
{c_5} = \frac{{L\sin \sigma }}{{m{v^2}{{\cos }^2}\theta \cos \psi }} - \frac{{h{\rho _0}{e^{ - \frac{h}{{{h_s}}}}}S({C_{L0}} + {C_{L\alpha }}\alpha )\sin \sigma }}{{2{h_s}m\cos \psi {{(\cos \theta )}^2}}}\\
{\rm{       }} + \theta (\frac{{{\rho _0}{e^{ - \frac{h}{{{h_s}}}}}S({C_{L0}} + {C_{L\alpha }}\alpha )\sin \sigma \sin \theta }}{{m\cos \psi {{(\cos \theta )}^3}}})\\
{\rm{       }} + \psi (\frac{{{\rho _0}{e^{ - \frac{h}{{{h_s}}}}}S({C_{L0}} + {C_{L\alpha }}\alpha )\sin \sigma \sin \psi }}{{2m{{(\cos \psi )}^2}{{(\cos \theta )}^2}}})\\
{\rm{        + }}\frac{{\sigma {\rho _0}{e^{ - \frac{h}{{{h_s}}}}}S({C_{L0}} + {C_{L\alpha }}\alpha )\cos \sigma }}{{2m\cos \psi {{(\cos \theta )}^2}}}\\
{\rm{       }} + \alpha (\frac{{{\rho _0}{e^{ - \frac{h}{{{h_s}}}}}S{C_{L\alpha }}\sin \sigma }}{{2m\cos \psi {{(\cos \theta )}^2}}})
\end{array}
\end{equation}
\begin{equation}
\begin{array}{l}
{c_6} = \frac{{\dot \sigma }}{{v\cos \theta \cos \psi }} - \frac{{\dot \sigma v}}{{{v^2}\cos \psi \cos \theta }}\\
{\rm{       }} - \frac{{\theta \dot \sigma \sin \theta }}{{v\cos \psi {{(\cos \theta )}^2}}} - \frac{{\psi \dot \sigma \sin \psi }}{{v{{(\cos \psi )}^2}\cos \theta }}\\
{\rm{       }} + \frac{{\dot \sigma }}{{v\cos \psi \cos \theta }}
\end{array}
\end{equation}


\bibliography{mybibfile}

\end{document}